\documentclass[12pt]{article}
\usepackage{times}    
\usepackage{lineno,amsthm, hyperref}
\usepackage{amsmath,amssymb,amscd,euscript}
\usepackage{dsfont}
\newtheorem{defi}{Definition}
\newtheorem{cor}{Corollary}
\newcommand{\dm}{{\mathcal{D}}}
\newcommand{\R}{{\mathbb{R}}}
\newcommand{\N}{{\mathbb{N}}}
\newcommand{\p}{{\mathbb{P}}}
\newcommand{\D}{{\mathbb{D}}}

\newcommand{\ep}{\epsilon}
\newcommand{\om}{\omega}

\newtheorem{lemma}{Lemma}[section]
\newtheorem{proposition}{Proposition}[section]
\newtheorem{theorem}{Theorem}[section]

\newtheorem{remark}{Remark}[section]

\allowdisplaybreaks \allowdisplaybreaks[4]

\begin{document}

\title{Malliavian differentiablity and smoothness of density for SDES with locally Lipschitz coefficients
\thanks{This work was supported by the Natural Sciences and Engineering Research Council of Canada under Grant DG-2018-04449.}
}


\author{Cristina Anton\\ 
Department of Mathematics and Statistics, MacEwan University, \\
103C, 10700-104 Ave., Edmonton, AB T5J 4S2, Canada \\
              Corresponding author, email: \texttt{popescuc@macewan.ca} 
}


\date{}

\maketitle
 
\begin{abstract}
 We study Malliavin differentiability for the solutions of a stochastic differential equation with drift of super-linear growth. Assuming we have a monotone drift with polynomial growth, we extend the work in \cite{ImkellerReisSalkeld:2019} and we prove Malliavin differentiability of any order. As a consequence of  this result, under the  H\"{o}rmander's hypothesis  we prove that the density of the solution's law with respect to the Lebesgue measure is infinitely differentiable. To avoid non-integrability problems due to the unbounded drift,  we follow an approach  based on the concepts of  Ray Absolute Continuity and Stochastic Gate{\^a}ux Differentiability. \\
\textbf{Keywords:} Malliavin  calculus, Monotone growth stochastic differential equation,  Wiener space
\end{abstract}
\section{Introduction}
We consider the stochastic differential equation (SDE)
\begin{equation}
dX(t)=b(X(t))dt+\sigma(X(t))dW(t),\quad
X(0)=x\in \mathbb{R}^d, t\in[0,T], T>0.\label{eq1}
\end{equation}
Here $W(t)$ is an $m-$ dimensional Brownian motion defined on the filtered complete probability space $\left(\Omega, \mathcal{F},\{\mathcal{F}_t\}_{t\ge 0}, \mathbb{P}\right)$, and $b:\mathbb{R}^d\rightarrow \mathbb{R}^d $,  $\sigma:\mathbb{R}^d\rightarrow \mathbb{R}^{d\times m}$. If $b$ and $\sigma$ are globally Lipschitz, $C^\infty$ and all their derivatives have polynomial growth, then the strong solution of \eqref{eq1} is Malliavin differentiable of any order \cite{KusuokaStroock:1985}. If in addition  the H\"{o}rmander's hypothesis holds at $x\in\mathbb{R}^d$, then the law of $X(t,0, x)$ is absolutely continuous with respect to the Lebesgue measure on $\mathbb{R}^d$ and it has an infinitely differentiable density \cite{KusuokaStroock:1985}. Here we obtain similar results under the following assumptions for the coefficients  $b$ and $\sigma$:
\begin{itemize}
\item[\bf{C:}] $b$  has bounded partial derivatives of any order $k\ge 2$ and $\sigma$ has bounded partial derivatives of any order $k\ge 1$.
\item[\bf{M:}] There exist $L>0$, such that for any $x_1, x_2\in \mathbb{R}^d $ we have
\begin{equation}
<x_1-x_2,b(x_1)-b(x_2)>\le L|x_1-x_2|^2\label{c2}
\end{equation}
\item[\bf{P:}] There exists $L_1\ge 0$ and $N\ge 1$ such that for any $x_1, x_2\in\mathbb{R}^d$ we have
\begin{equation}
|b(x_1)-b(x_2)|^2\le L_1(1+|x_1|^{2N-2}+|x_2|^{2N-2})|x_1-x_2|^2\label{pol}
\end{equation}
\end{itemize}
Thus, for $b$ we replace the globally Lipschitz conditions with  monotonicity ({\bf{M}}) and  polynomial growth ({ \bf{P}}) assumptions.

SDEs of this form are part of many stochastic models for real world applications, such as stochastic Ginzburg–Landau equation, stochastic Duffing–van der Pol oscillator, stochastic Lorenz equation,  stochastic Lotka–Volterra equations, and interest rate and volatility models in financial mathematics \cite{HutzenthalerJentzenKloeden:2011, KloedenPlaten:1992}. The convergence of some implicit numerical schemes for approximating the solutions of the SDE \eqref{eq1} was proven under monotonicity ({\bf{M}}) and  polynomial growth ({ \bf{P}}) assumptions \cite{TretyakovZhang:2013}. Our interest to study the properties of the density of the solution's law is motivated by potential applications to study the convergence rate of numerical schemes \cite{BallyTalay:1996}.

The global Lipschitz assumption is replaced in \cite{DeMarco:2011} by the assumptions that the coefficients are  smooth  and have bounded partial derivatives only on an open domain $D$. Estimations of the Fourier transform are used to show that the law of a strong solution of the
equation admits a smooth density on $D$, and under some additional conditions asymptotic estimates of the density are obtained. 

Using a characterization for Sobolev differentiability of random fields, in \cite{XieZhang:2016} the Sobolev regularity of strong solutions with respect to the initial value is established for SDEs with local Sobolev and super-linear growth coefficients. In \cite{RiedelScheutzow:2017} the existence of stochastic flows and semiflows  is discussed for SDEs driven by a rough path and with unbounded drift that grows at most linearly. 

First order Malliavin differentiability, parametric differentiability, and absolute continuity of the solution's law are proven in \cite{ImkellerReisSalkeld:2019} for SDEs with random coefficients with drifts satisfying locally Lipschitz  and monotonicity conditions. The concepts of Ray Absolute Continuity and Stochastic Gate{\^a}ux Differentiability are used to prove  Malliavin differentiability \cite{MastroliaPossamaiRveillac:2017}.   For SDE \eqref{eq1} under assumptions ({\bf{C}}), ({\bf{M}}), ({\bf{P}}) and the H\"{o}rmander's hypothesis, we combine this approach with techniques used in \cite{Nualart:2006}  and we prove Malliavin differentiability of any order and infinitely differentiability of the density of the solution's law.

The paper is organized as follows. In the next section we present some results of Malliavin calculus and elements of analysis on the Wiener space. Section \ref{mainsect} includes  preliminary results regarding Malliavin differentiability and the statement of the main result. In section \ref{sectlemas} we prove some technical lemmas. Based on these lemmas we prove Malliavin differentiability of any order for the solution of the SDE \eqref{eq1} in section \ref{sectproof}. Infinitely differentiability  of the density of the solution's law  is proven in section \ref{sect5}. 

\section{Preliminary results}\label{sect1}
Let $\left(\mathbb{R}^d, <\cdot,\cdot>, |\cdot|\right)$ be the d-dimensional real Euclidean space. If $A=(a_{ij})$ is a $d\times m$ matrix, we denote the Frobenius norm $|A|:=$trace$(AA^\top)^{1/2}=\left(\sum_{i=1}^{d}\sum_{j=1}^{m} a_{ij}^2\right)^{\frac{1}{2}}$. We denote by $\nabla f$ the gradient of a differentiable function $f:\mathbb{R}^d\rightarrow \mathbb{R}$.

We use the standard big $O$ and little $o$ notations for $f_n$, $f>0$: $f_n=O(f)$ means $\limsup_{n\rightarrow\infty}f_n/f\le C$, where $C>0$ is a constant independent of the limiting variable, and  $f_n=o(f)$ means $\lim_{n\rightarrow\infty}f_n/f=0$.

Let $p\ge 1$  and $T>0$. We inntroduce the following spaces.
\begin{itemize}
\item $C([0,T])$ be the Banach space of all uniformly continuous  and bounded functions $\phi:[0,T]\rightarrow \mathbb{R}$ endowed with the sup norm
$$
\|\phi\|_\infty=\sup_{x\in [0,T]}|\phi(x)|,
$$
and $\|\phi\|_{\infty,t}=\sup_{x\in [0,t]}|\phi(x)|$. Let $C_0([0,T])$ be the subspace of continuous functions that start at 0. 
\item Let $C_b^k(\mathbb{R}^d)$ the space of functions $\phi:\mathbb{R}^d\rightarrow \mathbb{R}$ consisting of all functions with bounded partial derivatives $\partial_x^i\phi(x)$, $0\le i\le k$ and with the norm
$$
\|\phi\|_k=\|\phi\|_\infty+\sum_{i=1}^k\sup_{x\in \mathbb{R}^d}|\partial_x^i\phi(x)|.
$$ 
Let $C_b^\infty(\mathbb{R}^d)=\cap_{k\ge 1} C_b^k(\mathbb{R}^d)$ and let $C_b^0(\mathbb{R}^d)$ be the space of continuous and bounded functions.
\item Let $C_p^\infty(\mathbb{R}^d)$ the set of all infinitely continuously differentiable functions $\phi:\mathbb{R}^d\rightarrow R$  such that $\phi$ and all it partial derivatives  have polynomial growth.
\item Let $L^p([0,T])=\{\phi:[0,T]\rightarrow \mathbb{R}, \|\phi\|_p=\left(\int_0^T |\phi (x)|^pdx\right)^{1/p}<\infty\}$. For any Hilbert space K with canonical norm $\|\cdot\|_K$ induced by the inner product and any $t\in[0,T]$ we set 
$$
L^p([t,T];K):=\left\{f:[t,T]\rightarrow K, \text{Borel measurable  and} \int_t^T\|f(s)\|_K^p ds<\infty\right\} 
$$
We denote $\mathcal{H}:=L^2([0,T];\mathbb{R}^m)$ and the canonical inner product is
$$
<f,g>_{\mathcal{H}}:=\int_0^T <f(s),g(s)>ds=\sum_{i=1}^m\int_0^Tf^i(s)g^i(s)ds,~f, g\in \mathcal{H}.
$$
\end{itemize}
Let  $\left(\Omega, \mathcal{F},\{\mathcal{F}_t\}_{t\ge 0}, \mathbb{P}\right)$ be a filtered probability space.
\begin{itemize}
\item  For $t\in [0,T]$ we define $L^p(\mathcal{F}_t,\mathbb{R}^d,\mathbb{P})=\{X:\Omega\rightarrow \mathbb{R}^d,$  $X$ is  $\mathcal{F}_t$ measurable and $\|X\|_p=E^\mathbb{P}[|X|^p]^{1/p}<\infty  \}$. For any separable Banach space $(E,\|\cdot\|)$, we denote $L^p(\Omega;E)=\{X:\Omega\rightarrow E$ , $X$ is  $\mathcal{F}$ measurable and $\|X\|_p=E^\mathbb{P}[\|X\|^p]^{1/p}<\infty  \}$.
Let $L^\infty$ be the subset of bounded random variables with norm $\|X\|_{L^\infty}=ess \sup_{\omega\in \Omega}|X(\omega)| $, and let $L^0(\mathcal{F}_t,\mathbb{R}^d,\mathbb{P})$ be the space of $X:\Omega\rightarrow \mathbb{R}^d,$  $X$ is  $\mathcal{F}_t$ measurable with the topology of convergence in probability.
\item Let $S^p([0,T],\mathbb{R}^d)=\{(Y_t)_{t\in [0,T]}$  stochastic processes, $Y_t\in\mathbb{R}^d$, that are $\{\mathcal{F}_t\}_{t\in [0,T]}$,   adapted and $\|Y\|_{S^p}=E^\mathbb{P}[|Y\|^p_\infty]^{1/p}=E^\mathbb{P}[\sup_{t\in [0,T]}|Y(t)|^p]^{1/p}$ $<\infty\}$. Let $S^\infty([0,T],\mathbb{R}^d)=\cap_{p\ge 1}S^p([0,T],\mathbb{R}^d)$.
\item We denote $H^p_{d,m}([0,T], \mathbb{R}^{d\times m}):=\biggl\{(Y_t)_{t\in [0,T]}$  stochastic processes, $Y_t\in\mathbb{R}^{d\times m}$, that are $\mathcal{F}$ predictable and $\|Y\|_{H^p_{d,m}}:=E\left[\left(\int_0^T\sum_{i=1}^d|Y_t^i|^2dt\right)^{p/2}\right]$ $<\infty.\biggl\}$. Put $H_m^p:=H^p_{1,m}$
\end{itemize}
\subsection{The Wiener space}
Let $\Omega=C_0\left([0,T],\mathbb{R}^m\right)=\{ \omega:[0,T]\rightarrow \mathbb{R}^m,\omega=(\omega_1,\ldots, \omega_m)^\top,$ $\omega_i\in C_0([0,T]), i=1,\ldots, m \} $ be the canonical Wiener space, and $W:=(W_t^1,\ldots, W_t^m)^\top_{t\in [0,T]}$ be the canonical Wiener process defined as $W_t^i(\omega):=\omega_t^i$ for any $\omega\in\Omega$, $i=1,\ldots, m$. We set $\mathcal{F}^0$ the natural filtration of $W$, $\mathbb{P}$ the Wiener measure, and $\mathcal{F}=\{\mathcal{F}_t\}_{t\in [0,T]}$ the usual augmentation (which is right-continuous and complete) of $\mathcal{F}^0$. In this setting $W$ is a standard Brownian motion. Let $H$ be the Cameron-Martin space:
\begin{align*}
H:=&\biggl\{h:[0,T]\rightarrow \mathbb{R}^m, h\in\Omega,\text{there exists } \dot{h}\in \mathcal{H}, \text{ such that } h(t)=\int_0^t \dot{h}(s)ds\text{ for any } t\in[0,T]\biggl\}
\end{align*}
For $h\in H$ we denote $\dot{h}$ a version of its Radon-Nykodym denisty with respect to the Lebesgue measure. H is a Hilbert space with the inner product $<h_1,h_2>_H:=<\dot{h}_1,\dot{h}_2>_{\mathcal{H}}$, for any $h_1,h_2\in H$, and with the associated norm $\|h\|_H^2=\|\dot{h}\|_{\mathcal{H}}^2$. For any Hilbert space $K$ we define $L^p(K)=\{f:\Omega\rightarrow K,$ $f$ is $\mathcal{F}_T$ measurable and $\|f\|^p_{L^p(K)}:=(E[\|f\|_K^p)^{1/p}<\infty\}$. Let
$$
W(h):=\int_0^T \dot{h}_s dW_s:=\sum_{i=1}^m \int_0^T \dot{h}_s^idW_s^i,~h\in H.
$$

For any $h\in H$ we define the shift operator $\tau_h:\Omega\rightarrow \Omega$, 
$
\tau_h(\omega):=\omega+h:=(\omega^1+h^1,\ldots, \omega^m+h^m)^\top.
$
We will list a few properties of the shift operator.
\begin{itemize}
\item For any $\mathcal{F}_T$ measurable $F:\Omega\rightarrow \mathbb{R}$, the mapping $h\rightarrow F\circ \tau_h$ is continuous in probability from $H$ to $L^0(\mathbb{R}^m):=\{f:\Omega\rightarrow \mathbb{R}, f$ is  $\mathcal{F}_T$ measurable$\}$ \cite[Lemma B2.1]{UstunelZakai:2000}. Consequently, $\tau_h$ is a continuous mapping on $\Omega$ for any $h\in H$.
\item If $X,Y:\Omega\rightarrow \mathbb{R}$ are $\mathcal{F}_T$ measurable and $X=Y$, $\mathbb{P}$ a.s., then for any $h\in H$ we have $X\circ \tau_h=Y\circ\tau_h$, $\mathbb{P}$ a.s.\cite[Appendix B.2]{UstunelZakai:2000}.
\item If $t\in [0,T]$ and $F:\Omega\rightarrow \mathbb{R}$ is $\mathcal{F}_t$ measurable then for any $h\in H$ it holds that $F\circ\tau_h=F\circ\tau_{h^t}$ $\mathbb{P}$ a.s., where
$$
\tau_{h^t}(s):=\left(\int_0^s \dot{h}^1(u)\mathds{1}_{[0,t]}(u)du,\ldots,\int_0^s \dot{h}^m(u)\mathds{1}_{[0,t]}(u)du \right)^\top.
$$
Consequently, $F\circ\tau_h$ is $\mathcal{F}_t$ measurable \cite[Lemma 3.3]{MastroliaPossamaiRveillac:2017}.
\item If $Y\in H^2_m$ and $h\in H$ then \cite[Lemma 3.4]{MastroliaPossamaiRveillac:2017}
\begin{equation}
\int_0^TY_sdW_s\circ\tau_h=\int_0^TY_s\circ\tau_hdW_s+\int_0^TY_s\circ\tau_h\dot{h}(s)ds, \mathbb{P} \text{ a.s.}.\label{helpto}
\end{equation}
\end{itemize}
We introduce the notation for the Dol\'eans-Dade exponential over $[0,T]$ of some sufficiently integrable $\mathbb{R}^m$ valued stochastic process $\{M(t)\}_{t\in [0,T]}$:
$$
\mathcal{E}(M)(t)=\exp\left(\int_0^t M(s)dW(s)-\frac{1}{2}\int_0^t|M(s)|^2ds\right)
$$
\begin{proposition}[ The Cameron-Martin Formula ]\cite[Appendix B.1]{UstunelZakai:2000}\\
Let $F:\Omega\rightarrow \mathbb{R}$ be $\mathcal{F}_t$ measurable and $h\in H$. When both sides are well defined 
$$
E[F\circ\tau_h]=E[F(\omega)\mathcal{E}(\dot{h})(T)]=E\left[F\exp\left(\int_0^T \dot{h}(s)dW(s)-\frac{1}{2}\int_0^T|\dot{h}(s)|^2ds\right)\right]
$$
Moreover, for any $h\in H$ and $p\ge 1$, $\mathcal{E}(\dot{h})(\cdot)\in S^p([0,T], \mathbb{R})$.
\label{cama}
\end{proposition}
\subsection{Results about Malliavin differentiability}
Following \cite{MastroliaPossamaiRveillac:2017} we set
\begin{align*}
\mathcal{S}:=&\biggl\{F:\Omega\rightarrow \mathbb{R}, F=f(W(h_1), \ldots,W(h_n)), f\in C_b^\infty(\mathbb{R}^n), \\
&h_i\in H, i=1,\ldots, n, \text{ for some } n\in\mathbb{N}, n\ge 1\biggl\} 
\end{align*}
For any $F\in \mathcal{S}$ we define the Malliavin derivative $\dm F:\Omega\rightarrow H$ by
$$
\dm F:=\sum_{i=1}^n\partial_i f(W(h_1), \ldots,W(h_n))h_i.
$$
We identify $\dm F$ with the stochastic process $\{\dm_tF\}_{t\in[0,T]}$, where $\dm_tF\in \mathbb{R}^m$ and
$$
\dm_tF(\omega)=\sum_{i=1}^n\partial_i f(W(h_1)(\omega), \ldots,W(h_n)(\omega))h_i(t), ~(t,\omega)\in [0,T]\times \Omega.
$$
$\mathcal{D}_t^jF$ will denote the $j$th component of $\mathcal{D}_tF$. 
We denote $\mathbb{D}^{1,p}$, $p\ge 1$ the closure of $\mathcal{S}$ with respect to the semi-norm
$$
\|F\|_{1,p}:=\left(E[|F|^p]+E[\|\dm F\|_H^p\right)^{1/p},
$$
and we set $\mathbb{D}^{1,\infty}=\cap_{p\ge 2}\mathbb{D}^{1,p}$.

The $k$th order Malliavin derivative $\mathcal{D}^kF:\Omega\rightarrow H^k$ is defined iteratively and its components are
$\mathcal{D}_{t_1,\ldots, t_k}^{j_1,\ldots, j_k}F:=\mathcal{D}_{t_k}^{j_k}\ldots \mathcal{D}_{t_1}^{j_1}F$,  with $(t_1,\ldots, t_k)^\top\in [0,T]^k$, $j_1, \ldots, j_k\in \{1, \ldots, m\}$. For the $N$th order Malliavin derivative, $\mathbb{{D}}^{N,p}$, $p\ge 1$ is the closure of $\mathcal {S}$ with the semi-norm
$$
\|F\|_{N,p}:=\left(E\left[|F|^p\right]+\sum_{i=1}^N E\left[\|\dm^i F\|_{H^i}^p\right]\right)^{1/p}=\left(\|F\|^p_{L^p(\Omega)}+\sum_{i=1}^N\|\mathcal{D}^iF\|^p_{L^p\left(\Omega;L^2([0,T]^i,\mathbb{R}^m)\right)}\right)^{\frac{1}{p}}.
$$
We set $\mathbb{{D}}^\infty=\cap_{p\ge 2}\cap_{i\ge 1}\mathbb{{D}}^{i,p}$

The definition of Malliavin derivative can be extended to mappings $G:\Omega\rightarrow E$, where $(E,\|\cdot\|_E)$ is a separable Banach space (\cite{ImkellerReisSalkeld:2019}). We consider the family
\begin{align*}
\mathcal{S}_E:=&\biggl\{G:\Omega\rightarrow E, G=\sum_{j=1}^kF_je_j, F_j\in S,~e_j\in E \text{ for some } k\in\mathbb{N}, k\ge 1\biggl\} 
\end{align*}
$\mathcal{S}_E$ is dense in $L^p(\mathcal{F};E;\mathbb{P})$ \cite{ImkellerReisSalkeld:2019}.
For any $G\in \mathcal{S}_E$ we define the Malliavin derivative $\dm G:\Omega\rightarrow H\otimes E$ by
$$
\dm G:=\sum_{j=1}^k\dm F_j\otimes e_j
$$
We denote $\mathbb{D}^{1,p}(E)$, $p\ge 1$ the closure of $\mathcal{S}_E$ with respect to the semi-norm
$$
\|G\|_{1,p,E}:=\left(E[\|G\|_E^p]+E[\|\|\dm G\|_H\|_E^p\right)^{1/p},
$$

Next, following \cite{MastroliaPossamaiRveillac:2017} we will present a different characterization of the spaces $\mathbb{D}^{1,p}$. We start with some definitions from \cite{ImkellerReisSalkeld:2019}. Let $E$ be a separable Banach space and $L(H,E)$ be the space of all bounded linear operators $V:H\rightarrow E$.
\begin{defi}
A measurable map $f:\Omega\rightarrow E$ is said to be Ray Absolutely Continuous if for any $h\in H$ there exists a measurable mapping $\tilde{f}_h:\Omega\rightarrow E$ such that $\tilde{f}_h(\omega)=f(\omega)$, $\mathbb{P}$ a.s., and that for any $\omega\in \Omega$, $t\rightarrow \tilde{f}_h(\omega+th) $ is absolutely continuous on any compact subset of $\mathbb{R}$.
\end{defi}
\begin{defi}
A measurable map $f:\Omega\rightarrow E$ is said to be Stochastically Gate{\^a}ux differentiable if there exists a measurable mapping $F:\Omega\rightarrow L(H,E)$ such that for any $h\in H$,
$$
\frac{f(\omega+\epsilon h)-f(\omega)}{\epsilon}\overset{\mathbb{P}}{\rightarrow}F(\omega)[h]~\text{as }\epsilon\rightarrow 0.
$$
\end{defi}
\begin{defi}
Let $p>1$. A measurable map $f\in L^p(\Omega;E)$ is said to be Strong Stochastically Gate{\^a}ux differentiable if there exists a measurable mapping $F:\Omega\rightarrow L(H,E)$ such that for any $h\in H$,
$$
\underset{\epsilon\rightarrow 0}{\lim}E\left[\left\|\frac{f(\omega+\epsilon h)-f(\omega)}{\epsilon}-F(\omega)[h]\right\|\right]=0.
$$
\end{defi}
\begin{theorem}\cite[Theorem 3.10]{ImkellerReisSalkeld:2019}\\
\label{thg1}
Let $p>1$. The space $\mathbb{D}^{1,p}(E)$ is equivalent with the space of all random variables $f:\Omega\rightarrow E$ such that $f\in L^p(\Omega;E)$ is Ray Absolutely Continuous, Stochastically Gate{\^a}ux differentiable, and the Stochastic Gate{\^a}ux derivative $F:\Omega\rightarrow L(H,E)$ is $F\in L^p(\Omega;L(H,E)) $.
\end{theorem}
\begin{theorem}\cite[Theorem 3.13]{ImkellerReisSalkeld:2019}\\
\label{th3.13}
Let $p>1$. The space $\mathbb{D}^{1,p}(E)$ is equivalent with the space of all random variables $f:\Omega\rightarrow E$ such that $f\in L^p(\Omega;E)$ is Strong Stochastically Gate{\^a}ux differentiable and have measurable mappings  $F\in L^p(\Omega;L(H,E)) $.
\end{theorem}
\begin{cor}\cite[Corollary 3.14]{ImkellerReisSalkeld:2019}\\
\label{cori}
If a measurable map $f:\Omega\rightarrow E$ is Stochastically Gate{\^a}ux differentiable and for   $\delta>0$
$$
\underset{\epsilon\in (0,1]}{\sup}E\left[\left\|\frac{f(\omega+\epsilon h)-f(\omega)}{\epsilon}\right\|^{1+\delta}\right]<\infty,
$$
then $f$ is Malliavin differentiable (and so is Ray Absolutely Continuous)
\end{cor}
\begin{lemma}
\label{lemaA21}
Let $q\ge 2$ and $f\in L^q([0,T],\mathbb{D}^{1,q})$. Then for any $p\in(1,q)$ and for any $h\in H$
\begin{equation}
\underset{\epsilon\rightarrow 0}{\lim}E\left[\int_0^T\left|\frac{f_s\circ\tau_{\epsilon h}-f_s}{\epsilon}-<\dm f_s,h>_{H}\right|^p ds\right]=0.\label{la21}
\end{equation}
\end{lemma}
The proof is included in appendix A.
\begin{remark}
\label{rem1}
Lemma \ref{lemaA21} can be easily extended to $f:\Omega\times[0,T]\rightarrow \R ^d$ such that for  each component $f^i\in L^q\left([0,T],\mathbb{D}^{1,q}\right)$, $i=1,\ldots, d$.  We have for any $h\in H$ and any $p\in (1,q)$:
\begin{equation}
\underset{\epsilon\rightarrow 0}{\lim}E\left[\int_0^T\left|\frac{f_s\circ\tau_{\epsilon h}-f_s}{\epsilon}-\dm f_s[h]\right|^p ds\right]=0,\label{la21bis}
\end{equation}
where
$$
\dm f_s[h]=\int_0^T \dm_t f_s\dot{h}(t)dt=(<\dm f_s^1,h>_{H}, \ldots,<\dm f_s^d,h>_{H}>)^\top\in\R^d 
$$
\end{remark}
\begin{lemma}
\label{lemaA22}
Let $q\ge 2 $ and $f:[0,T]\times \Omega\rightarrow \R$ such that
\begin{equation}
E\left[\left(\underset{t\in[0,T]}{\sup}|f_t|\right)^q\right]<\infty,~E\left[\left(\underset{t\in[0,T]}{\sup}\int_0^T|\dm_s f_t|^2ds\right)^{q/2}\right]<\infty.\label{eql1}
\end{equation}
Then we have
\begin{align}
&\underset{t\in[0,T]}{\sup}\left|\frac{f_t\circ\tau_{\epsilon h}-f_t}{\epsilon}-<\dm f_t,h>_H\right|~\overset{\p}{\underset{\epsilon\rightarrow 0}{\rightarrow}}~0\label{probi11}\\
&\underset{\epsilon\in(0,1]}{\sup}E\left[\underset{t\in[0,T]}{\sup}\left|\frac{f_t\circ\tau_{\epsilon h}-f_t}{\epsilon}\right|^2\right]<\infty\label{probi21}
\end{align}
\end{lemma}
The proof is included in appendix B.
\begin{remark}
\label{rem2}
Lemma \ref{lemaA22} can be easily extended to and $f:[0,T]\times \Omega\rightarrow \R^d$ such that
\begin{equation}
E\left[\left(\underset{t\in[0,T]}{\sup}|f_t|\right)^q\right]<\infty,~E\left[\left(\underset{t\in[0,T]}{\sup}\int_0^T|\dm_s f_t|^2ds\right)^{q/2}\right]<\infty.\label{eql}
\end{equation}
Then we have
\begin{align}
&\underset{t\in[0,T]}{\sup}\left|\frac{f_t\circ\tau_{\epsilon h}-f_t}{\epsilon}-\dm f_t[h]\right|~\overset{\p}{\underset{\epsilon\rightarrow 0}{\rightarrow}}~0\label{probi1}\\
&\underset{\ep\in(0,1]}{\sup}E\left[\underset{t\in[0,T]}{\sup}\left|\frac{f_t\circ\tau_{\epsilon h}-f_t}{\epsilon}\right|^2\right]<\infty\label{probi2}
\end{align}
\end{remark}
\section{The main result}\label{mainsect}
From assumption {\bf{C}} we get that for any $R>0$ there exists $K_R>0$ such that
\begin{equation}
|b(x_1)-b(x_2)|\le K_R|x_1-x_2|\text{ for any } |x_1|\le R,  |x_2|\le R .\label{c1}
\end{equation}
Furthermore, assumption {\bf{C}} implies that  $\sigma$ and  $\nabla_xb$ are globally Lipschitz,  and assumption  {\bf{M}} implies that for any $x,y\in\R ^d$
\begin{equation}
y^\top\nabla_xb(x)y\le L|y|^2.\label{c4}
\end{equation}
From assumptions {\bf{C}} and {\bf{M}},  the Cauchy-Schwartz  and Young inequalities we get for any $\epsilon>0$ and any  $x\in\R ^d$
\begin{align}
&2<x,b(x)>\le 2L|x|^2+2|x||b(0)|\le(2L+\epsilon)|x|^2+\frac{|b(0)|^2}{\epsilon},\label{cc1}\\ 
&|\sigma(x)|^2\le L_2(1+\epsilon)|x|^2+|\sigma(0)|^2\left(1+\frac{1}{\epsilon}\right)\le C(1+|x|^{2}).\label{us4}
\end{align}

From assumption {\bf{C}}, (\ref{cc1}), (\ref{us4}) and Theorems 3.6 in \cite{Mao:2011} we know that there exists a unique global solution $X(t,0,x)$ of the SDE (\ref{eq1}).
Also the solution $\left(X(t,0,x)\right)_{t\ge 0}$ is $\{\mathcal{F}_t\}_{t\ge 0}$ adapted and we have
\begin{equation}
E\left(\int_0^T |X(t,0,x)|^2dt\right)<\infty.\label{mom1}
\end{equation} 
Moreover from Theorems 9.1 and 9.5 in \cite{Mao:2011} we know that $\left(X(t,0,x)\right)_{t\ge 0}$ is a time homogeneous Markov process, and from Theorem 4.1 in \cite{Mao:2011} we know that  for any $p\ge  2$ there exists a constant $\alpha_p>0$ such that we have
\begin{equation}
E\left[|X(t,0,x)|^p\right]\le 2^{\frac{p-2}{2}}\left(1+|x|^p\right)e^{p\alpha_pt}, \quad t\in [0,T].\label{mom2}
\end{equation}

From Theorem 2.2 in \cite{ImkellerReisSalkeld:2019}, for any $p\ge 2$ we have $X\in S^p:=S^p([0,T],\mathbb{R}^d)$ and there exists $C>0$ depending on $p$, $b$ and $\sigma$ such that:
\begin{equation}
E\left[\sup_{t\in[0,T]}|X(t,0,x)|^p \right]< C(|x|^p+1).\label{mom3}
\end{equation}

Notice that using assumptions {\bf{C}} and {\bf{P}} and proceeding as for \eqref{us4} we can show that there exist $C>0$ such that we have for any $x\in\mathbb{R}^d$
\begin{align}
|\nabla_x b(x)|^2&\le C(1+|x|^2)\label{us2}\\
|b(x)|^2&\le C(1+|x|^{2N})\label{us3}
\end{align}
From these inequalities, \eqref{mom3} and assumption {\bf{C}} , for any multi-index $\alpha=(\alpha_1,\ldots,\alpha_d)$  and any $p\ge 2$ we get
\begin{align}
&E\left[\sup_{t\in[0,T]}|\partial_\alpha b(X(t))|^p \right]< \infty,~
E\left[\sup_{t\in[0,T]}|\partial_\alpha \sigma(X(t))|^p \right]< \infty.\label{mom31}
\end{align}

From Corollary 3.5 in \cite{ImkellerReisSalkeld:2019} we know that $X$ is Malliavin differentiable and $\mathcal{D}_sX(t)=0$ for $T\ge s>t\ge 0$ while for $0\le s\le t\le T$
\begin{equation}
\dm_sX(t)=\sigma(X(s))+\int_s^t\nabla_xb(X(r))\dm_sX(r)dr+\int_s^t\nabla_x\sigma(X(r))\dm_sX(r)dW(r)\label{mad1}
\end{equation} 
Notice that from Theorem 3.21 in \cite{ImkellerReisSalkeld:2019} we know that SDE \eqref{mad1} has a unique solution in  $S^p\left([0,T],L^2([0,T])\right)$, for any $p\ge 2$. Thus we have $\dm X\in S^p\left([0,T],L^2([0,T])\right)$:
\begin{equation}
E\left[\left(\sup_{t\in[0,T]}\int_0^T|\dm_sX(t)|^2ds\right)^{p/2} \right]< \infty\label{malmom1}
\end{equation}
Thus from \eqref{mom3} and \eqref{malmom1}  we have $ X\in \mathbb{D}^{1,p}(S^p)$ for any $p\ge 2$, so $ X\in \mathbb{D}^{1,\infty}(S^p)=\cap_{p\ge 2}\mathbb{D}^{1,p}(S^p)$. Thus, similarly with the case in Theorem 2.2.1 in \cite{Nualart:2006} of globally Lipschitz coefficients, for any $t\in[0,T]$ we have $ X(t)\in \mathbb{D}^{1,\infty}$.

For any fixed $s\in [0,T]$ and $i=1,\ldots, m$, from \eqref{mad1} we know that $\dm_s^iX$ verifies the linear SDE
\begin{equation}
\dm_s^iX(t)=\sigma^i(X(s))+\int_s^t\nabla_xb(X(r))\dm_s^iX(r)dr+\int_s^t\nabla_x\sigma(X(r))\dm_s^iX(r)dW(r),
\end{equation} 
so from \eqref{us4}, \eqref{c4}, \eqref{mom1}, and Theorem 2.5 in \cite{ImkellerReisSalkeld:2019} we have for any $p\ge 2$
\begin{align*}
&E\left[\underset{t\in [0,T]}{\sup} |\dm_s^iX(t)|^p\right]\le C E\left[|\sigma^i(X(s,0,x))|^p\right]\le C \left(1+E\left[|X(s))|^{p}\right]\right)
\end{align*}
Thus from \eqref{mom2} we get for any $p\ge 2$
\begin{align}
&\underset{s\in [0,T]}{\sup}E\left[\underset{t\in [0,T]}{\sup} |\dm_sX(t, 0,x)|^p\right]\le c_{1,p}(T)(1+|x|^p)<\infty\label{sav1}, 
\end{align}
where $c_{1,p}(T)>0$ depends on $p$, $T$, $b$ and $\sigma$.
This and \eqref{mom2} implies that for any $t\in [0,T]$ and any $p\ge 2$
\begin{align}
&\|X(t, 0,x)\|_{1,p}^p=E[|X(t, 0,x)|^p]+E\left[\left|\int_0^T|\dm_sX(t, 0,x)|^2 ds\right|^{p/2}\right]\le C_{1,p}(T)(1+|x|^p), \label{mom3-1}
\end{align}
where $C_{1,p}(T)>0$ depends on $p$, $T$, $b$ and $\sigma$.
Here we obtain the following extension of this result.
\begin{theorem}
\label{mainh}
Let $\{X(t,0,x)\}_{t\in[0,T]}$ with $X(0)=x$ be the solution of SDE \eqref{eq1} under assumptions  {\bf{C}}, {\bf{M}}, and {\bf{P}}. Then $X^i(t)$ belongs to $\D ^\infty$ for all $t\in[0,T]$, and $i=1,\ldots, d$. Moreover, for any $t\in [0,T]$, $p\ge 2$, $k=1,2,\ldots, $ there exist $C_{k,p}(T), \beta_{k,p}>0$ depending on $p$, $k$, $T$, $b$, $\sigma$ such that
\begin{align}
\|X(t, 0,x)\|_{k,p}^p\le C_{k,p}(T)(1+|x|^{\beta_{k,p}})\label{dolfif}
\end{align}. 
\end{theorem}
 The proof of this theorem will be given in section \ref{sectproof}.
\section{Some technical lemmas}
\label{sectlemas}
Before we prouve Theorem \ref{mainh} we need a few technical results.
Let $\{X(t)\}_{t\in [0,T]}$ with  $X(0)=x$ be the solution of SDE \eqref{eq1} under assumptions  {\bf{C}}, {\bf{M}}, and {\bf{P}}.
Let $r\in [0,T)$ and $\{Y(t)\}_{r\le t\le T}\subset \R^d$ be the solution of the linear SDE
\begin{align}
&Y(t)=\alpha(r,\omega)+\int_r^t U(u,\omega)du+\sum_{i=1}^m\int_r^tV^i(u,\omega)dW^i(u)\notag\\
&+\int_r^t\nabla_x b(X(u))Y(u)du+\sum_{i=1}^m\int_r^t\nabla_x\sigma^i(X(u))Y(u)dW^i(u), ~t\ge r,\label{eq3}
\end{align}
where $\{\alpha (r)\}_{r\in [0,T]}\subset \R^d$, $\{U(t)\}_{t\in [0,T]}\subset \R^d$ and $\{V(t)\}_{t\in [0,T]}\subset \R^{d\times m}$   are adapted and continuous processes. We assume that for any $r\in[0,T]$, $\alpha (r)\in \D ^{1,\infty}(\R ^d)$ and for any $t\in[r,T]$, $U(t)\in \D ^{1,\infty}(\R ^d)$, $V(t)\in \D ^{1,\infty}(\R ^{d\times m})$. Moreover, for any $p\ge 2$ we assume that there exists the positive constants $C$, $p_{\alpha_1}$, $p_{\alpha_2}$, $p_{U_1}$, $p_{U_2}$, $p_{V_1}$, $p_{V_2}$ depending on $p$, $T$, $b$ and $\sigma$ such that we have 
\begin{align}
&E\left[\underset{r\in [0,T]}{\sup}|\alpha (r)|^p\right]\le C(1+|x|^{p_{\alpha_1}})<\infty,\label{alp1}\\
&\underset{s\in [0,T]}{\sup}E\left[\underset{r\in [s,T]}{\sup}|D_s\alpha (r)|^p\right]\le C(1+|x|^{p_{\alpha_2}})<\infty,\label{alp3}\\
&E\left[\underset{s\in [r,T]}{\sup}|U(s)|^p\right]\le C(1+|x|^{p_{U_1}})<\infty,\label{U1}\\
&\underset{s\in [0,T]}{\sup}E\left[\underset{u\in [r\vee s,T]}{\sup}|D_s U(u)|^p\right]\le C(1+|x|^{p_{U_2}})<\infty,\label{U3}\\
&E\left[\left(\underset{s\in [r,T]}{\sup}|V(s)|^p\right)\right]\le C(1+|x|^{p_{V_1}})<\infty,\label{V1}\\
&\underset{s\in [0,T]}{\sup}E\left[\underset{u\in [r\vee s,T]}{\sup}|D_s V(u)|^p\right]\le C(1+|x|^{p_{V_2}})<\infty.\label{V3}
\end{align}
From Theorem 2.5 in \cite{ImkellerReisSalkeld:2019} we know that there exists a unique solution of equation \eqref{eq3} in $S^p\left([r,T],\R ^d\right)$,  $p\ge 2$. Moreover, from the proof of Theorem 2.5 in \cite{ImkellerReisSalkeld:2019} and  H\"{o}lder inequality  we have for any $p\ge 2$
\begin{align*}
&E\left[\sup_{t\in[r,T]}|Y(t)|^p \right]\le C\biggl(E\left[|\alpha (r)|^p\right]++E\left[\left(\int_r^T|U(s)|\right)^p\right]+E\left[\left(\int_r^T|V(s)|^2\right)^{p/2}\right]\biggl)\\
&\le C\biggl(E\left[|\alpha (r)|^p\right]+E\left[\underset{s\in [r,T]}{\sup}|U(s)|^p\right]+E\left[\underset{s\in [r,T]}{\sup}|V(s)|^p\right]\biggl)
\end{align*}
Hence \eqref{alp1}, \eqref{U1},\eqref{V1} imply  that for any $p\ge 2$ there exists the positive constants $C$, $p_{Y_1}$ depending on $p$, $T$, $b$ and $\sigma$ such that we have
\begin{equation}
E\left[\sup_{t\in[r,T]}|Y(t)|^p \right]\le C(1+|x|^{p_{Y_1}})< \infty.\label{mom3y}
\end{equation}.
\begin{remark}
Notice that for any $p\ge 2$, \eqref{alp3}, \eqref{U3}, \eqref{V3} and  H\"{o}lder inequality imply that there exist the positive constants  $C$, $C_1$, $p_{\alpha_3}$, $p_{U_3}$, $p_{V_3}$ depending on $p$, $T$, $b$ and $\sigma$ such that we have  
\begin{align}
&E\left[\underset{r\in [0,T]}{\sup}\left(\int_0^T|D_s\alpha (r)|^2ds\right)^{p/2}\right]\le C(1+|x|^{p_{\alpha_3}})<\infty,\label{alp2}\\
&E\left[\underset{u\in [r,T]}{\sup}\left(\int_0^T|D_s U(u)|ds\right)^{p}\right]\le C_1E\left[\underset{u\in [r,T]}{\sup}\left(\int_0^T|D_s U(u)|^2ds\right)^{p/2}\right]\notag\\
&\le C(1+|x|^{p_{U_2}})<\infty,\label{U2}\\
&E\left[\underset{u\in [r,T]}{\sup}\left(\int_0^T|D_s V(u)|^2ds\right)^{p/2}\right]\le C(1+|x|^{p_{V_3}})<\infty.\label{V2}
\end{align}
\end{remark}
We have assumed that  $\alpha^i (r)\in \D^{1,\infty}(\R)$, $i=1,\ldots, d$. From Theorem 4.1 and Lemma 4.2 in \cite{MastroliaPossamaiRveillac:2017}  we get for any $p\ge 1$:
 \begin{align*}
 &\underset{\epsilon\rightarrow 0}{\lim}E\left[\left|\frac{\alpha^i (r)\circ\tau_{\epsilon h}-\alpha^i (r)}{\epsilon}- <\dm \alpha^i (r),h>_H\right|^p\right]=0
 \end{align*}
As in Remark \ref{rem1}  we can show that this imply
\begin{align}
& \underset{\epsilon\rightarrow 0}{\lim}E\left[\left|\frac{\alpha (r)\circ\tau_{\epsilon h}-\alpha (r)}{\epsilon}- \dm \alpha (r) [h]\right|^p\right]=0\label{con1a}
\end{align}
 where
 \begin{align*}
&\dm \alpha (r)[h]=\int_0^T \dm_t \alpha (r)\dot{h}(t)dt=\int_0^r \dm_t \alpha (r)\dot{h}(t)dt=(<\dm \alpha^1 (r),h>_{H}, \ldots,<\dm \alpha^d (r),h>_{H}>)^\top
\end{align*}

Similarly, since for any $t\in [r,T]$, $U(t)\in \D^{1,\infty}(\R^d)$, $V(t)\in \D ^{1,\infty}(\R ^{d\times m})$ from \eqref{U1}, \eqref{U2}, \eqref{V1}, \eqref{V2}, and Remarks \ref{rem1} and \ref{rem2} we have for any $h\in H$ and any $p>1$:
\begin{align}
&\underset{\epsilon\rightarrow 0}{\lim}E\left[\int_r^T\left|\frac{U(s)\circ\tau_{\epsilon h}-U(s)}{\epsilon}-\dm U(s)[h]\right|^p ds\right]=0,\label{con1U}\\
&\underset{\epsilon\in(0,1]}{\sup}E\left[\underset{s\in[r,T]}{\sup}\left|\frac{U(s)\circ\tau_{\epsilon h}-U(s)}{\epsilon}\right|^2\right]<\infty,\label{con3U}\\
&\underset{\epsilon\rightarrow 0}{\lim}E\left[\int_r^T\left|\frac{V^j(s)\circ\tau_{\epsilon h}-V^j(s)}{\epsilon}-\dm V^j(s)[h]\right|^p ds\right]=0,\label{con1V0}\\
&\underset{\epsilon\in(0,1]}{\sup}E\left[\underset{s\in[r,T]}{\sup}\left|\frac{V^j(s)\circ\tau_{\epsilon h}-V^j(s)}{\epsilon}\right|^2\right]<\infty,~ j=1,\ldots, m,\label{con3V0}
\end{align}
where
\begin{align*}
&\dm U(s)[h]=\int_0^s \dm_t U(s)\dot{h}(t)dt=(<\dm U^1(s),h>_{H}, \ldots,<\dm U^d(s),h>_{H}>)^\top\\
 &\dm V^j(s)[h]=\int_0^s \dm_t V^j(s)\dot{h}(t)dt=(<\dm V^{1,j}(s),h>_{H}, \ldots,<\dm V^{d,j}(s),h>_{H}>)^\top,
\end{align*}
and $ \dm V(s)[h]=\int_0^T \dm_t V(s)\dot{h}(t)dt$ is a $d\times m $ matrix with elements $ \dm V(s)[h]^{i,j}=<\dm V^{i,j}(s),h>_{H}$.

Let $B(t, \omega):=\nabla_x b(X(t)(\omega))$. From assumption {\bf{C}}, $X(t)\in \D ^{1,\infty}$, and chain rule we get   $B^{i,j}(t)\in \D ^{1,\infty}$ and $\dm _s B^{i,j}(t)=\sum_{k=1}^d\partial_{j,k} b^i(X(t))\dm _s X^k(t)$, for any $t\in [r,T]$ and any $i,j=1,\ldots, d$. Using assumption {\bf{C}}, \eqref{mom3}, \eqref{malmom1}, and Remarks \ref{rem1} and \ref{rem2} we have for any $h\in H$ and any $p>1$: 
\begin{align}
&\underset{\epsilon\rightarrow 0}{\lim}E\left[\int_r^T\left|\frac{B(s)\circ\tau_{\epsilon h}-B(s)}{\epsilon}-\dm B(s)[h]\right|^p ds\right]=0,\label{con1B}\\
&\underset{\epsilon\in(0,1]}{\sup}E\left[\underset{s\in[r,T]}{\sup}\left|\frac{B(s)\circ\tau_{\epsilon h}-B(s)}{\epsilon}\right|^2\right]<\infty,\label{con3B}
\end{align}
where $\dm B(s)[h]=\int_0^T \dm_t B(s)\dot{h}(t)dt=\int_0^s \dm_t B(s)\dot{h}(t)dt$ is a $d\times d $ matrix with elements $ \dm B(s)[h]^{i,j}=<\dm B^{i,j}(s),h>_{H}$. For any $p\ge 2$ from \eqref{malmom1} and \eqref{mom3y} we have
\begin{align}
&E\left[\left(\int_r^T\left(\int_0^T|\dm_ s \nabla_x b(X(u))Y(u)|^2ds\right)^{1/2}du\right)^p\right]\notag\\
&\le C \left(E\left[\left(\underset{u\in [r,T]}{\sup}|Y(u)|^{2p}\right)\right]\right)^{1/2}\left(E\left[\left(\underset{u\in [r,T]}{\sup}\int_0^T|\dm_ s X(u)|^2ds\right)^{p}\right]\right)^{1/2}<\infty\label{con4B1}
\end{align}
Moreover, since
\begin{align*}
&0\le E\left[\left(\int_r^T\left|\frac{\left(B(u,\omega+\epsilon h)-B(u,\omega)\right)Y(u)}{\epsilon}-\int_0^T\dm_s B(u,\omega) \dot{h}(s) dsY(u)\right|du\right)^2\right]\\
&\le C \left( E\left[\int_r^T\left|\frac{B(u,\omega+\epsilon h)-B(u,\omega)Y(u)}{\epsilon}-\int_0^T\dm_s B(u,\omega) \dot{h}(s) ds\right|^4du\right]\right)^{1/2}\left( E\left[\underset{u\in [r,T]}{\sup}|Y(u)|^4\right]\right)^{1/2},
\end{align*}
putting $p=4$ in \eqref{con1B} and \eqref{mom3y} we get
\begin{align}
&\underset{\epsilon\rightarrow 0}{\lim}E\left[\left(\int_r^T\left|\frac{\left(\nabla_x b(X(u)(\omega+\epsilon h))-\nabla_x b(X(u)(\omega)\right))Y(u)}{\epsilon}-\int_0^u\dm_s \nabla_x b(X(u)) \dot{h}(s)Y(u) ds\right|du\right)^2\right]\notag\\
&= 0\label{con5B}
\end{align}
Thus  from this and \eqref{mom3y} with p=4 we get
\begin{align}
&E\left[\left(\int_r^T\left|\frac{\left(B(u,\om +\ep h))-B(u,\om))\right)Y(u)(\om )}{\ep} \right| du\right)^2\right]=o(1)+2d^2T\|\dot{h}\|_2^2\notag\\
&\sum_{i,j=1}^d \left(E\left[\underset{u\in[r,T]}{\sup}|Y(u)(\om )|^4du\right]\right)^{1/2}\left(E\left[\int_r^T\left(\int_0^T|\dm_s B^{i,j}(u)(\om )|^2ds\right)^2du \right]\right)^{1/2}\label{con6B}
\end{align}

Similarly, let $\Sigma(t, \omega):=\nabla_x \sigma (X(t)(\omega))$. From assumption {\bf{C}}, $X(t)\in \D ^{1,\infty}$, and chain rule we get  $\Sigma^{i,j,k}(t)\in \D ^{1,\infty}$ and $\dm _s \Sigma^{i,j,k}(t)=\sum_{l=1}^d\partial_{j,l} \sigma^{i,k}(X(t))\dm _s X^l(t)$, for  for any $t\in [r,T]$ and any $i,j=1,\ldots, d$, $k=1,\ldots, m$. Using assumption {\bf{C}}, \eqref{malmom1},   and Remarks \ref{rem1} and \ref{rem2} we have for any $h\in H$ and any $p>1$:
\begin{align}
&\underset{\epsilon\rightarrow 0}{\lim}E\left[\int_r^T\left|\frac{\Sigma^k(s)\circ\tau_{\epsilon h}-\Sigma^k(s)}{\epsilon}-\dm \Sigma^k(s)[h]\right|^p ds\right]=0,\label{con1s0}\\
&\underset{\epsilon\in(0,1]}{\sup}E\left[\underset{s\in[r,T]}{\sup}\left|\frac{\Sigma^k(s)\circ\tau_{\epsilon h}-\Sigma^k(s)}{\epsilon}\right|^2\right]<\infty,\label{con3s0}
\end{align}
where $\dm \Sigma^k(s)[h]=\int_0^T \dm_t \Sigma^k(s)\dot{h}(t)dt=\int_0^s \dm_t \Sigma^k(s)\dot{h}(t)dt$ is a matrix in $\R^{d\times d}$ with elements $ \dm \Sigma^k(s)[h]^{i,j}=<\dm \Sigma^{i,j,k}(s),h>_{H}$.
Using assumption {\bf{C}}, \eqref{malmom1}, and \eqref{mom3y}, for any $p\ge 2$ we have 
\begin{align}
&E\left[\left(\int_r^T\int_0^T|\dm_ s \nabla_x \sigma^k(X(u))Y(u)|^2ds du\right)^{p/2}\right]<\infty\label{con4s}
\end{align}
Analogously with \eqref{con5B}, putting $p=4$ in \eqref{con1s0} and \eqref{mom3y} we get for any $k=1,\ldots, m$
\begin{align}
&E\left[\int_r^T\left|\frac{\left(\nabla_x \sigma^k(X(u)(\omega+\epsilon h))-\nabla_x \sigma^k(X(u)(\omega)\right))Y(u)}{\epsilon}-\int_0^T\dm_s \nabla_x \sigma^k(X(u)) \dot{h}(s)dsY(u) \right|^2du\right]\notag\\
&\rightarrow 0\text{ as } \epsilon\rightarrow 0\label{con5s0}
\end{align}
From this and \eqref{mom3y} with $p=4$  we get
\begin{align}
&E\left[\int_r^T\left|\frac{\left(\Sigma^i(u,\om +\ep h)-\Sigma^i(u,\om)\right)Y(u)(\om)}{\ep}\right|^2du\right]\notag\\
&=o(1)+2\|\dot{h}\|_2^2\sum_{j,k=1}^d\left(E\left[\int_r^T\left(\int_0^T\left|\dm_s \Sigma^{i,j,k}(u)\right|^2ds\right)^2du\right]\right)^{1/2}\left(E\left[\underset{u\in[r,T]}{\sup}|Y(u)|^4du\right]\right)^{1/2}\label{con6s0}
\end{align}
\begin{lemma}
\label{lemh1}
Let $p\ge 2$, $r\in[0,T)$, and $Y$ be the solution of the SDE \eqref{eq3} under assumptions  {\bf{C}}, {\bf{M}}, {\bf{P}},  and \eqref{alp1}-\eqref{V3}. For $(s,t)\in [0,T]^2$, let $(M_s(t))$ be defined by the $d\times m$ matrix of $L^2([r,T])$- valued SDEs
\begin{align}
&M_s(t)=\dm_s\alpha (r,\om)\mathds{1}_{\{r\ge s\}}+ \left(\nabla_x\sigma(X(s))Y(s)+V(s)\right)\mathds{1}_{\{r\le s\}}+\int_{s\vee r}^t \dm _s U(u,\omega)du\notag\\
&+\sum_{i=1}^m\int_{s\vee r}^t\dm _s V^i(u,\omega)dW^i(u)+\int_{s\vee r}^t\dm_s \nabla_x b(X(u))Y(u)du+\int_{s\vee r}^t\nabla_x b(X(u))M_s(u)du\notag\\
&+\sum_{i=1}^m\int_{s\vee r}^t\dm _s \nabla_x\sigma^i(X(u))Y(u)dW^i(u)+\sum_{i=1}^m\int_{s\vee r}^t\nabla_x\sigma^i(X(u))M_s(u)dW^i(u),\label{eq4bis}
\end{align}
for $s\le t$, and $M_s(t)=0$ for $s>t$. Then a unique solution exists in $S^p\left([r,T],L^2([r,T])\right)$ for \eqref{eq4bis}. Moreover we have
\begin{equation}
E\left[\left(\underset{t\in [r\vee s,T]}{\sup}\int_0^T|M_s(t)|^2ds\right)^{p/2}\right]<\infty.
\end{equation} 
\end{lemma}
The proof is included in appendix C.
\begin{lemma}
\label{lemah2}
Let $Y$ be the solution of SDE \eqref{eq3} under assumptions  {\bf{C}}, {\bf{M}}, {\bf{P}},  and \eqref{alp1}-\eqref{V3}. We have 
\begin{equation}
\underset{\ep \in (0,1)}{\sup}E\left[\underset{t\in  [r,T]}{\sup}\left|P_{\ep}(t)(\om)\right|^2\right]<\infty,~ P_{\ep}(t)(\om):=\frac{Y(t)(\om+\ep h)-Y(t)(\om)}{\ep}
\end{equation} 
\end{lemma}
The proof is included in appendix D.
\begin{lemma}
\label{lemh3}
Let $Y$ be the solution of SDE \eqref{eq3}, $M_s$ be the solution of \eqref{eq4bis}  under assumptions  {\bf{C}}, {\bf{M}}, {\bf{P}}, and \eqref{alp1}-\eqref{V3}, and let $h\in H$. We have 
\begin{equation}
\underset{\ep \rightarrow 0}{\lim}\left\|\frac{Y(\cdot)(\om+\ep h)-Y(\cdot)(\om) }{\ep}-M^h(\cdot)(\om)\right\|_\infty\overset{ \p }{\rightarrow}0,~M^h(t)(\om):=\int_0^t M_s(t)(\om)\dot{h}(s)ds
\end{equation} 
i.e. $Y$ is Stochastically Gate{\^a}ux differentialble.
\end{lemma}
The proof is included in appendix E.
\begin{theorem}
\label{thsg}
Let $Y$ be the solution of SDE \eqref{eq3}, $M_s$ be the solution of \eqref{eq4bis}  under assumptions  {\bf{C}}, {\bf{M}}, {\bf{P}},   and \eqref{alp1}-\eqref{V3}, and let $h\in H$. We have 
\begin{equation}
\underset{\ep \rightarrow 0}{\lim}E\left[\left\|\frac{Y(\cdot)(\om+\ep h)-Y(\cdot)(\om) }{\ep}-\int_0^\cdot M_s(\cdot) \dot{h}(s)ds\right\|_\infty\right]=0
\end{equation} 
i.e. $Y$ is Strong Stochastically Gate{\^a}ux differentialble.
\end{theorem}
\begin{proof}
From Lemma \ref{lemh3} we have convergence in probability. From Lemma \ref{lemh1} for any $p\ge 2$ we get
\begin{align*}
&E\left[\|M^h(\om)\|_\infty^p\right]:=E\left[\underset{t\in [r,T]}{\sup}|M^h(t)(\om)|^p\right]\le E\left[\left(\underset{t\in [r,T]}{\sup}|M^h(t)(\om)|\right)^p\right]\\
&\le E\left[\left(\underset{t\in [r,T]}{\sup}\left( \int_0^t | M_s(t)(\om)|^2ds\right)^{1/2}\left(\int_0^t|\dot{h}(s)|^2ds\right)^{1/2}\right)^p\right]\\
&\le \|\dot{h}\|_2^{p} E\left[\left(\underset{t\in [r,T]}{\sup} \int_0^T| M_s(t)(\om)|^2ds\right)^{p/2}\right]<\infty
\end{align*}
Combining  with Lemma \ref{lemah2} we  have 
\begin{equation*}
\underset{\ep \in (0,1)}{\sup}E\left[\left\|\frac{Y(t)(\om+\ep h)-Y(t)(\om)}{\ep}\right\|_\infty^2\right]<\infty, ~E\left[\|M^h(\om)\|_\infty^2\right]<\infty,
\end{equation*} 
so we apply Corollary \ref{cori} to conclude.
\end{proof}
 \begin{theorem}
\label{lemah1}
Let $\{X(t)\}_{t\in [0,T]}$ with  $X(0)=x$ be the solution of SDE \eqref{eq1}.  Let $r\in [0,T)$ and $\{Y(t)\}_{r\le t\le T}\subset \R^d$ be the solution of the linear SDE \eqref{eq3}.
Then under assumptions  {\bf{C}}, {\bf{M}}, {\bf{P}} and \eqref{alp1}-\eqref{V3}, $Y\in\D^{1,p}(S^p)$, for any $p\ge 2$, and the Malliavin derivative satisfies for $r\vee s\le t\le T$ and $j=1,\ldots, m$
\begin{align}
&\dm _s^j Y(t)=\dm_s^j\alpha (r)\mathds{1}_{r\ge s}+ \left(\nabla_x\sigma^j(X(s))Y(s)+V^j(s)\right)\mathds{1}_{r\le s}+\int_{r\vee s}^t \dm _s^j U(u,\omega)du\notag\\
&+\sum_{i=1}^m\int_{r\vee s}^t\dm _s^j V^i(u,\omega)dW^i(u)+\int_{r\vee s}^t\dm_s^j\nabla_x b(X(u))Y(u)du+\int_{r\vee s}^t\nabla_x b(X(u))\dm_s^j Y(u)du\notag\\
&+\sum_{i=1}^m\int_{r\vee s}^t\dm _s^j\nabla_x\sigma^i(X(u))Y(u)dW^i(u)+\sum_{i=1}^m\int_{r\vee s}^t\nabla_x\sigma^i(X(u))\dm_s^j Y(u)dW^i(u),\label{eq4}
\end{align}
and $\dm _s Y(t)=0$ for $s\vee r>t$.
Moreover  for any $p\ge 2$ there exist the positive constants $C$, $p_{Y_2}$  depending on $p$, $T$, $b$ and $\sigma$ such that we have we have
\begin{align}
&E\left[\left(\underset{t\in [r,T]}{\sup}\int_0^T|\dm _s^jY(t)|^2ds\right)^{p/2}\right]<\infty,\label{usold}\\
&\underset{s\in [0,T]}{\sup}E\left[\underset{t\in [r\vee s,T]}{\sup}|\dm _s^jY(t)|^p\right]\le C(1+|x|^{p_{Y_2}})<\infty.\label{sav2}
\end{align} 
\end{theorem}
 \begin{proof}
 In Theorem \ref{thsg}  we have shown that $Y$ is Strong Stochastically Gate{\^a}ux differentiable and $\dm ^h Y(\cdot):=M^h(\cdot)$ has finite moments of any order $p\ge 2$, so  $Y\in\D^{1,p}(S^p)$.  
The Malliavin derivatives satisfies equation \eqref{eq4}, and from Lemma \ref{lemh1} we know that the equation \eqref{eq4} has a unique solution.
 
We get inequality \eqref{usold} from Lemma \ref{lemh1}.  To prove \eqref{sav2}, for any fixed $r, s\in [0,T]$ and $j=1,\ldots, m$ we apply Theorem 2.5 in \cite{ImkellerReisSalkeld:2019} to the linear SDE \eqref{eq4} with  $r\vee s\le t\le T$  and we have for any $p\ge 2$ from assumption {\bf{C}},  and Jensen inequality 
\begin{align*}
&E\left[\underset{t\in [r\vee s,T]}{\sup}|\dm _s^jY(t)|^p\right]\le C\biggl( E\left[\left(|\dm_s^j\alpha (r)\mathds{1}_{r\ge s}+ \left(\nabla_x\sigma^j(X(s))Y(s)+V^j(s)\right)\mathds{1}_{r\le s}|\right)^p\right]\\
&+E\left[\left(\int_{r\vee s}^T |\dm _s^j U(u)+\dm_s^j\nabla_x b(X(u))Y(u)|du\right)^p\right]\\
&+E\left[\left(\int_{r\vee s}^T |\dm _s^j V(u)+\dm _s^j\nabla_x\sigma(X(u))Y(u)|^2du\right)^{p/2}\right]\biggl)\\
&\le C\biggl(E\left[\underset{r\in [s,T]}{\sup}|\dm_s^j\alpha (r)|^p\right]+E\left[|Y(s)|^p\right]\mathds{1}_{r\le s}+E\left[|V(s)|^p\right]\mathds{1}_{r\le s}\\
&+E\left[\underset{u\in [r\vee s,T]}{\sup} |\dm _s^j U(u)|^p\right]+E\left[\underset{u\in [r\vee s,T]}{\sup}|\dm_s^j X(u)|^{2p}\right]+E\left[\underset{u\in [r\vee s,T]}{\sup}|Y(u)|^{2p}\right]\\
&+E\left[\underset{u\in [r\vee s,T]}{\sup}  |\dm _s^j V(u)|^{p}\right]+E\left[\underset{u\in [r\vee s,T]}{\sup}|\dm _s^j X(u)|^{2p}\right]+E\left[\underset{u\in [r\vee s,T]}{\sup}|Y(u)|^{2p}\right]\biggl)
 \end{align*}
 Hence
 \begin{align*}
 &\underset{s\in [0,T]}{\sup}E\left[\underset{t\in [r\vee s,T]}{\sup}|\dm _s^jY(t)|^p\right]\le\notag\\
 &\le C\biggl(\underset{s\in [0,T]}{\sup}E\left[\underset{r\in [s,T]}{\sup}|\dm_s^j\alpha (r)|^p\right]+E\left[\underset{s\in [r,T]}{\sup}|Y(s)|^p\right]+E\left[\underset{s\in [r,T]}{\sup}|V(s)|^p\right]\notag\\
&+\underset{s\in [0,T]}{\sup}E\left[\underset{u\in [r\vee s,T]}{\sup} |\dm _s^j U(u)|^p\right]+\underset{s\in [0,T]}{\sup}E\left[\underset{u\in [r\vee s,T]}{\sup}|\dm_s^j X(u)|^{2p}\right]+E\left[\underset{u\in [r,T]}{\sup}|Y(u)|^p\right]\notag\\
&+\underset{s\in [0,T]}{\sup}E\left[\underset{u\in [r\vee s,T]}{\sup}  |\dm _s^j V(u)|^{p}\right]+\underset{s\in [0,T]}{\sup}E\left[\underset{u\in [r\vee s,T]}{\sup}|\dm _s^j X(u)|^{2p}\right]+E\left[\underset{u\in [r,T]}{\sup}|Y(u)|^{2p}\right]\biggl)
 \end{align*}
 Thus from \eqref{sav1}, \eqref{alp3}, \eqref{U3}, \eqref{V1}, \eqref{V3}, \eqref{mom3y} we get \eqref{sav2}.
 \end{proof}
\section{Proof of Theorem \ref{mainh} }
\label{sectproof}
\begin{proof}
 We know that for any $t\in[0,T]$, $j=1,\ldots, m$, $i=1,\ldots, d$, we have $X^i(t)\in \D^{1,\infty}$ and $\dm _r^j X(t)=0$  for $T\ge r>t\ge 0$ while for $0\le r\le t\le T$ verifies the following linear stochastic differential equation (see \eqref{mad1}):
\begin{equation*}
\dm_r^jX(t)=\sigma^j(X(r))+\int_r^t\nabla_xb(X(u))\dm_r^jX(u)du+\int_r^t\nabla_x\sigma(X(u))\dm _r^j X(u)dW(u)
\end{equation*} 
We will recursively apply Theorem \ref{lemah1} to this linear equation. Notice that the previous equation is similar with \eqref{eq3} with $\alpha (r,\omega)=\sigma^j(X(r)(\om))$ and $U(u,\om)=0$, $V(u,\om)=0$, $u\in[0,T]$. For any $p\ge 2$ we have from \eqref{mom3} and \eqref{us4}:
\begin{align*}
&E\left[\underset{r\in [0,T]}{\sup}|\sigma^j(X(r))|^p\right]\le CE\left[\underset{r\in [0,T]}{\sup}(1+|X(r)|^2)^{p/2}\right]\\
&\le  C(1+E\left[\underset{r\in [0,T]}{\sup}|X(r)|^{p}\right]\le C(1+|x|^p)<\infty
\end{align*}
From \eqref{sav1} and assumption {\bf{C}} we get for any $p\ge 2$
\begin{align}
&\underset{s\in [0,T]}{\sup}E\left[\underset{r\in [s,T]}{\sup}\left|\dm _s \sigma^j(X(r)) \right|^p\right]\le\underset{s\in [0,T]}{\sup}E\left[\underset{r\in [s,T]}{\sup}|\nabla_x\sigma^j(X(r))|^p| \dm _s X(r) |^p\right]\notag\\
&\le C \underset{s\in [0,T]}{\sup}E\left[\underset{r\in [s,T]}{\sup} |\dm _s X(r) |^p\right]\le C(1+|x|^p)<\infty.\label{sav4}
\end{align}
Thus \eqref{alp1}-\eqref{V3} hold and from Theorem \ref{lemah1} we get that $\dm_{rs}^{jk}X(t):=\dm_s^k\dm_r^jX(t)=0$  for $s\vee r>t$ and for  $s\vee r\le t\le T$ it verifies the equation
\begin{align}
&\dm _{rs}^{jk} X(t)=\dm_s^k\sigma^j(X(r)(\om))\mathds{1}_{r\ge s}+ \nabla_x\sigma^k(X(s))\dm_r^jX(s)\mathds{1}_{r\le s}\notag\\
&+\int_{r\vee s}^T\dm_s^k\nabla_x b(X(u))\dm_r^jX(u)du+\sum_{i=1}^m\int_{r\vee s}^T\dm _s^k\nabla_x\sigma^i(X(u))\dm_r^jX(u)dW^i(u)\notag\\
&+\int_{r\vee s}^T\nabla_x b(X(u))\dm_{rs}^{jk}X(u)du+\sum_{i=1}^m\int_{r\vee s}^T\nabla_x\sigma^i(X(u))\dm_{rs}^{jk}X(u)dW^i(u)\label{sav5}
\end{align}
For any $p\ge 2$, $r\in[0,T]$, $j,k=1,\ldots,m$, we also have from \eqref{usold} 
\begin{align}
&E\left[\left(\underset{t\in [r,T]}{\sup}\int_0^T|\dm _{rs}^{jk} X(t)|^2ds\right)^{p/2}\right]<\infty.\label{sav6}
\end{align} 
Moreover \eqref{sav2} implies that for any $p\ge 2$ there exist $C, p_2>0$ such that
\begin{align}
&\underset{s\in [0,T]}{\sup}E\left[\underset{t\in [r\vee s,T]}{\sup}|\dm _{rs}^{jk}X(t)|^p\right]\le C(1+|x|^{p_{2}})\label{sav7}
\end{align}
This and H\"{o}lder inequality imply that for any $t\in [0,T]$, $p\ge 2$ we have
\begin{align*}
&E\left[\left(\int_0^T\int_0^T |\dm _{rs}^{jk}X(t)|^2drds\right)^{p/2}\right]\le C \int_0^T\int_0^TE\left[ |\dm _{rs}^{jk}X(t)|^p\right]drds\\
&\le C\underset{r,s\in [0,T]}{\sup} E\left[|\dm _{rs}^{jk}X(t)|^p\right]\le C(1+|x|^{p_{2}})<\infty.
\end{align*}
Thus $X(t) \in\D^{2,\infty}$ for any $t\in [0,T]$, and using also \eqref{mom3} and \eqref{mom3-1} we get that for any $p\ge 2$ there exist $C_{2,p}(T), \beta_{2,p}>0$ depending on $p$, $T$, $b$, $\sigma$ such that
\begin{align}
\|X(t)\|_{2,p}^p\le C_{2,p}(T)(1+|x|^{\beta_{2,p}}).\label{dolfi2}
\end{align}

Next, the linear equation \eqref{sav5} is similar with  \eqref{eq3} with $\alpha (s,\omega)=\dm_s^k\sigma^j(X(r)(\om))\mathds{1}_{r\ge s}+ \nabla_x\sigma^k(X(s))\dm_r^jX(s)\mathds{1}_{r\le s}$, $U(u)=\dm_s^k\nabla_x b(X(u))\dm_r^jX(u)$, and $V(u)=\dm _s^k\nabla_x\sigma(X(u))\dm_r^jX(u)$, $u\in[0,T]$. From Assumption C, \eqref{sav1}, \eqref{sav4} we get that for any $p\ge 2$ there exist $C,p_2>0$ such that we have
\begin{align*}
&E\left[\underset{r\vee s\in [0,T]}{\sup}|\dm_s^k\sigma^j(X(r))\mathds{1}_{r\ge s}+ \nabla_x\sigma^k(X(s))\dm_r^jX(s)\mathds{1}_{r\le s}|^p\right]\le C(1+|x|^{p_3})<\infty.
\end{align*}
From Assumption {\bf{C}}, Young inequality, \eqref{sav1}, and \eqref{sav7} we get that for any $p\ge 2$ there exist $C,p_2>0$ such that we have
\begin{align*}
&\underset{u\in [0,T]}{\sup}E\left[\underset{s\vee r\in [u,T]}{\sup}|\dm_u\alpha (s)|^p\right]\le C\underset{u\in [0,T]}{\sup}E\left[\underset{s\vee r\in [u,T]}{\sup}|\dm_u\left(\dm_s^k\sigma^j(X(r))\right)\mathds{1}_{r\ge s}|^p\right]\\
&+ C\underset{u\in [0,T]}{\sup}E\left[\underset{s\vee r\in [u,T]}{\sup}|\dm_u\left( \nabla_x\sigma^k(X(s))\dm_r^jX(s)\right)\mathds{1}_{r\le s}|^p\right]\\
&\le C\underset{u\in [0,T]}{\sup}E\left[\underset{ r\in [u,T]}{\sup}|D^2\sigma^j(X(r))\dm_uX(r)\dm_s^k X(r)|^p\mathds{1}_{r\ge s}\right]\\
&+C\underset{u\in [0,T]}{\sup}E\left[\underset{ r\in [u,T]}{\sup}|\nabla_x\sigma^j(X(r))\dm_u\dm_s^k X(r)|^p\mathds{1}_{r\ge s}\right]\\
&+C\underset{u\in [0,T]}{\sup}E\left[\underset{s\in [u,T]}{\sup}| D^2\sigma^k(X(s))\dm_uX(s)\dm_r^jX(s)|^p\mathds{1}_{r\le s}\right]\\
&+C\underset{u\in [0,T]}{\sup}E\left[\underset{s\in [u,T]}{\sup}| \nabla_x\sigma^k(X(s))\dm_u\dm_r^jX(s)|^p\mathds{1}_{r\le s}\right]\\
&\le C\left(\underset{u\in [0,T]}{\sup}E\left[\underset{ r\in [u,T]}{\sup}|\dm_uX(r)|^{2p}\right]+\underset{s\in [0,T]}{\sup}E\left[\underset{ r\in [s,T]}{\sup}|\dm_s^k X(r)|^{2p}\right]\right)\\
&+C\underset{u\in [0,T]}{\sup}E\left[\underset{ r\in [u\vee s,T]}{\sup}|\dm_u\dm_s^k X(r)|^p\right]\\\
&+C\left(\underset{u\in [0,T]}{\sup}E\left[\underset{s\in [u,T]}{\sup}| \dm_uX(s)|^{2p}\right]+\underset{r\in [0,T]}{\sup}E\left[\underset{s\in [r,T]}{\sup}\dm_r^jX(s)|^{2p}\right]\right)\\
&+C\underset{u\in [0,T]}{\sup}E\left[\underset{s\in [u\vee r,T]}{\sup}| \dm_u\dm_r^jX(s)|^p\right]\le C(1+|x|^{p_3})<\infty
\end{align*}

For $U(u)=\dm_s^k\nabla_x b(X(u))\dm_r^jX(u)$, from assumption {\bf{C}}, Young inequality and \eqref{sav1} we have for any $p\ge 2$ 
\begin{align*}
&E\left[\underset{u\in [r\vee s,T]}{\sup}|U(u)|^p\right]=E\left[\underset{u\in [r\vee s,T]}{\sup}|D^2 b(X(u))\dm_s^kX(u)\dm_r^jX(u)|^p\right]\\
&\le CE\left[\underset{u\in [r\vee s,T]}{\sup}|\dm_s^kX(u)\dm_r^jX(u)|^p\right]\le C(1+|x|^{p_4})<\infty,
\end{align*}
where $C,p_4>0$ depend on $p$, $T$, $b$ and $\sigma$. Moreover, from  assumption {\bf{C}}, Young inequality, \eqref{sav1}, and \eqref{sav7} we get that there exist $C,p_4>0$ depending on $p$, $T$, $b$ and $\sigma$ such that we have
\begin{align*}
&\underset{v\in [0,T]}{\sup}E\left[\underset{u\in [r\vee s\vee v,T]}{\sup}|\dm_v U(u)|^p\right]=\underset{v\in [0,T]}{\sup}E\left[\underset{u\in [r\vee s\vee v,T]}{\sup}|\dm_v(D^2 b(X(u))\dm_s^k X(u)\dm_r^jX(u))|^p\right]\\
&\le C\underset{v\in [0,T]}{\sup}E\left[\underset{u\in [r\vee s\vee v,T]}{\sup}|D^3 b(X(u))|^p|\dm_v X(u)|^p|\dm_s^k X(u)|^p|\dm_r^jX(u)|^p\right]\\
&+C\underset{v\in [0,T]}{\sup}E\left[\underset{u\in [r\vee s\vee v,T]}{\sup}|D^2 b(X(u))|^p| \dm_v\dm_s^k X(u)|^p\dm_r^jX(u)|^p\right]\\
&+C\underset{v\in [0,T]}{\sup}E\left[\underset{u\in [r\vee s\vee v,T]}{\sup}|D^2 b(X(u))|^p| \dm_s^k X(u)|^p|\dm_v\dm_r^jX(u)|^p\right]\le C(1+|x|^{p_4})<\infty.
\end{align*}

Finally, proceeding in the same way for $V(u)=\dm _s^k\nabla_x\sigma(X(u))\dm_r^jX(u)$ we can show that for any $p\ge 2$ there exist $C,p_5>0$ depending on $p$, $T$, $b$ and $\sigma$ such that we have
\begin{align*}
&E\left[\left(\underset{u\in [r\vee s,T]}{\sup}|V(u)|^p\right)\right]\le C(1+|x|^{p_5})<\infty,\\
&\underset{v\in [0,T]}{\sup}E\left[\underset{u\in [r\vee s\vee v,T]}{\sup}|\dm_v V(u)|^p\right]\le C(1+|x|^{p_5})<\infty.
\end{align*}
Thus we can apply again Theorem \ref{lemah1} to the linear equation \eqref{sav5} and show that $X(t) \in\D^{3,\infty}$, for any $t\in [0,T]$.
As in the proof of Theorem 2.2.2 in \cite{Nualart:2006} we can show by induction on $N$ that $X(t)\in\D^{N,\infty}$ for any $N\in \N^*$, $t\in [0,T]$. Thus  $X(t)\in \D ^\infty$ for any $t\in [0,T]$. 
\end{proof}

\section{Applications}
\label{sect5}
From Theorem 4.9 in \cite{ImkellerReisSalkeld:2019} we know that under assumptions {\bf{C}}, {\bf{M}}, {\bf{P}} the matrix valued SDE
\begin{equation}
J(t)=I_d+\int_0^t\nabla_x b\left(X(s,0,x)\right)J(s)ds+\int_0^t\nabla_x\sigma\left(X(s,0,x)\right)J(s)dW(s), \label{jaco}
\end{equation}
$t\in[0,T]$, has a unique solution $J\in S^p\left([0,T],\mathbb{R}^{d\times d}\right)$, $p\ge 2$, and for any $t\in[0,T]$ the map $x\rightarrow X(t,0,x)$ is differentiable $\mathbb{P}$ a.s. and as $\epsilon\rightarrow 0$
\begin{equation}
\frac{X(t, 0,x+\epsilon h)(\omega)-X(t,0,x)(\omega)}{\epsilon}\rightarrow hJ(t)(\omega)~\mathbb{P}~ a.s. . 
\end{equation}

To show that the Jacobian matrix $J(t)$ is $\mathbb{P}$ a.s. invertible for any choice of $t\in [0,T]$, we consider an extra assumption.
\begin{itemize}
\item[\bf{J:}]  There exist $L_3>0$ such that for any $y\in\mathbb{R}^d$ we have
\begin{equation}
y^\top \nabla_x by>-L_3|y|^2.
\end{equation}
\end{itemize}
From Theorem 2.5 and Proposition 4.13 in \cite{ImkellerReisSalkeld:2019} we know that under assumptions {\bf{C}}, {\bf{M}}, {\bf{P}}, {\bf{J}} the matrix valued SDE
\begin{align}
&K(t)=I_d-\int_0^tK(s)\left[\nabla_x b\left(X(s,0,x)\right)-<\nabla_x\sigma,\nabla_x\sigma>\left(X(s,0,x)\right)\right]ds\notag\\
&-\int_0^tK(s)\nabla_x\sigma\left(X(s,0,x\right)dW(s),\quad t\in[0,T], \label{jacin}
\end{align}
has a unique solution $K\in S^p\left([0,T],\mathbb{R}^{d\times d}\right)$, $p\ge 2$, and we have $K(t)J(t)=I_d$ for all $t\in[0,T]$ $\mathbb{P}$ a.s.. Consequently, the Jacobian matrix $J(t)$ is $\mathbb{P}$ a.s. invertible for any choice of $t\in[0,T]$, and $J(t)^{-1}=K(t)$ $\mathbb{P}$ a.s..

Let $J_s(t)=J(t)J(s)^{-1}$, $t>s$. Under  assumptions {\bf{C}}, {\bf{M}}, {\bf{P}}, {\bf{J}} from Proposition 5.1 in \cite{ImkellerReisSalkeld:2019} we know that we have
\begin{equation}
J_s(t)=I_d+\int_s^t\nabla_x b\left(X(r,0,x)\right)J_s(r)dr+\int_s^t\nabla_x\sigma\left(X(r,0,x)\right)J_s(r)dW(r), \label{jacofun}
\end{equation}
and the Malliavin derivative of $X$ can be expressed for $t>s$ as $D_sX(t)=J_s(t)\sigma(X(s))$. $J_s(t)$ is the fundamental matrix of the linear SDE \eqref{mad1}.
The Malliavin matrix $Q(t)$ is defined by
\begin{align}
Q(t,x)&:=\int_0^tD_sX(t,0,x)D_sX(t,0,x)^\top ds=J(t)C(t,x)J(t)^\top\label{matQ}\\
C(t,x)&:=\int_0^t J(s)^{-1}\sigma(X(s,0,x))\sigma(X(s,0,x))^\top(J(s)^{-1})^\top ds\label{matc}
\end{align}

The Lie bracket of the $C^1(\mathbb{R}^d,\mathbb{R}^d)$ vector fields $V=\sum_{i=1}^d V^i \frac{\partial}{\partial x_i}$, $U=\sum_{i=1}^d U^i \frac{\partial}{\partial x_i}$ is defined as
$[V,U](x)=\partial U(x)V(x)-\partial V(x) U(x)$, where $\partial U=(\partial_iU^j)_{i,j=1,\ldots d}$, $\partial V=(\partial_iV^j)_{i,j=1,\ldots d}$ are the Jacobian matrices of $U$ and $V$ respectively. Let us denote $\sigma^0=b-\frac{1}{2}\sum_{i=1}^{m}\sum_{j=1}^{d}\sigma_j^{i}\partial_j \sigma^i$ and let $\sigma^0$, $\ldots$, $\sigma^m$ be the corresponding vector fields 
$$
\sigma^0(x)=\sum_{i=1}^d \sigma_0^i(x)\frac{\partial}{\partial x_i},\quad \sigma^j(x)=\sum_{i=1}^d \sigma^j_{i}(x)\frac{\partial}{\partial x_i}, \quad j=1,\ldots, m.
$$
We construct by recurrence the sets $\Sigma_0=\{\sigma^j, j=1,\ldots,m \}$, $\Sigma_k=\{[\sigma^j, V],  j=0,\ldots,m, V\in\Sigma_{k-1}\},~ k\ge 1$, $\Sigma_\infty=\cup_{k=1}^\infty \Sigma_k$.
We denote by $\Sigma_k(x)$ the subset of $\mathbb{R}^m$ obtained by freezing the variable $x\in \mathbb{R}^d$ in the vector fields of $\Sigma_k$. Note than $Span\{\Sigma_0(x)\}=\mathbb{R}^m$ is equivalent with $\sigma\sigma^\top(x)>0$ which is the ellipticity assumption in $x$. Under this assumption  and assumptions {\bf{C}}, {\bf{M}}, and {\bf{P}},  from Theorem 5.2 in \cite{ImkellerReisSalkeld:2019} we know that the law of $X(t,0, x)$ is absolutely continuous with respect to the Lebesgue measure on $\mathbb{R}^d$.

For $x\in\mathbb{R}^d$ we consider the H\"{o}rmander's hypothesis which is weaker than the ellipticity  assumption:
\begin{itemize}
\item[\bf{H(x): }] The vector space $Span\{\Sigma_\infty(x)\}=\mathbb{R}^d$ 
\end{itemize}
\begin{theorem}
 Let $\{X(t,0,x)\}_{t\in [0,T]}$ with  $X(0)=x$ be the solution of SDE \eqref{eq1}.  Suppose that  {\bf{H(x)}}, {\bf{C}}, {\bf{M}}, {\bf{P}}, {\bf{J}}  hold. Then for any $t\in [0,T]$ the law of $X(t,0, x)$ is absolutely continuous with respect to the Lebesgue measure on $\mathbb{R}^d$ and has an infinitely differentiable density.
\end{theorem}
\begin{proof}
From \eqref{mom3-1} we know that $X(t, 0,x)\in \D ^{1,p}$ for any $p\ge 2$. According to \cite[Theorem 2.2.1]{Nualart:2006}, to prove that  the law of $X(t,0, x)$ is absolutely continuous with respect to the Lebesgue measure, it is enough to show that  the Malliavin matrix $Q(t,x)$ is invertible a.s. Based on assumption {\bf{H(x)}}, \eqref{mom31} and the formulas \eqref{matQ}, \eqref{matc} for the Malliavin matrix $Q(t,x)$,  the proof can be done as the proof of Theorem 2.3.2 in  \cite{Nualart:2006}. Moreover,  using \eqref{dolfif} in Theorem \ref{mainh} and proceeding as in the proof of Theorem 2.3.3 in \cite{Nualart:2006}, we can show the density with respect to the Lebesgue measure is infinitely differentiable.
\end{proof}


\appendix
\section*{Appendix A: The proof of lemma \ref{lemaA21}}
\label{apendixa}
{\small
\begin{proof}
Since $f\in L^q([0,T],\mathbb{D}^{1,q})$, for all $s\in [0,T]$, we have $f_s:\Omega\rightarrow \mathbb{R}$, $f_s\in \mathbb{D}^{1,q}$. From Theorem \ref{thg1} we  get that $f_s$ is Ray Absolutely Continuous and Stochastically Gate{\^a}ux differentiable, so for any $h\in H$ and any $p\in(1,q)$, Theorem 4.1 in \cite{MastroliaPossamaiRveillac:2017} yields
\begin{equation}
{\underset{\epsilon\rightarrow 0}{\lim}}E\left[\left|\frac{f_s\circ\tau_{\epsilon h}-f_s}{\epsilon}-<\dm f_s,h>_{H}\right|^p\right]=0.\label{eqsa}
\end{equation}
Moreover, from equation (4.2) in \cite{MastroliaPossamaiRveillac:2017} for any $\epsilon\in (0,1)$ and any $h\in H$  we have
\begin{equation}
\frac{f_s\circ\tau_{\epsilon h}-f_s}{\epsilon}=\frac{1}{\epsilon}\int_0^\epsilon <\dm f_s\circ\tau_{u h} ,h>_{H}du,~ \mathbb{P} \text{ a.s.}\label{eqsh}
\end{equation}
From $f\in L^q([0,T],\mathbb{D}^{1,q})$ we get $\|f_s\|_{1,q}<\infty$ and this implies $f\in \mathbb{D}^{1,q}(E)$, where $E=L^q([0,T])$.

From Theorem \ref{thg1} we  get that $f$ is Ray Absolutely Continuous and Stochastically Gate{\^a}ux differentiable. Thus $f\in L^q(\Omega, E)$ and for any $p\in(1,q)$ using H\"{o}lder inequality and we can show that $f\in L^p(\Omega, L^p[0,T])$. Moreover, for any $\epsilon>0$ there exists $\delta(\epsilon)>0$ such that we have
\begin{align*}
\mathbb{P}\left(\int_0^T\left|\frac{f_s\circ\tau_{\epsilon h}-f_s}{\epsilon}-F_s[h]\right|^q ds>\delta(\epsilon)\right)<\epsilon,
\end{align*}
where $F:\Omega\rightarrow L(H,E)$, $F\in L^q(\Omega, L(H,E))$.
Using this and H\"{o}lder inequality, it is easy to show that for any $p \in(1, q)$, any $\epsilon\in(0,1)$,  and any $h\in H$, we have  $\frac{f\circ\tau_{\epsilon h}-f}{\epsilon}\in L^p(\Omega, L^p[0,T])$ and 
\begin{align}
\frac{f\circ\tau_{\epsilon h}-f}{\epsilon}\overset{\mathbb{P}}{\underset{\epsilon\rightarrow 0}{\rightarrow}}F[h]\label{pr1}
\end{align}

Next we show that  for any $p\in (1,q)$ the family $\left\{\left\|\frac{f\circ\tau_{\epsilon h}-f}{\epsilon}\right\|_{L^p[0,T]}^p \right\}_{\epsilon\in (0,1)}=\left\{\int_0^T\left|\frac{f_s\circ\tau_{\epsilon h}-f_s}{\epsilon}\right|^p ds\right\}_{\epsilon\in (0,1)}$ is uniformly integrable in $L^p(\Omega, L^p[0,T])$ .
We fix any $p\in (1,q)$, $h\in H$ and $\eta>0$ such that $r:=p+\eta\in (p,q)$.
From Jensen inequality we get
\begin{align}
&E\left[\left(\int_0^T\left|\frac{f_s\circ\tau_{\epsilon h}-f_s}{\epsilon}\right|^pds\right)^{(p+\eta)/p}\right]\le 
 E\left[\int_0^T\left|\frac{f_s\circ\tau_{\epsilon h}-f_s}{\epsilon}\right|^rds\right]\label{equi1}
\end{align}
Proceeding as in the proof of Lemma 4.2 in \cite{MastroliaPossamaiRveillac:2017}, using \eqref{eqsh},  H\"{o}lder  inequality, the Cameron-Martin Formula   we obtain
\begin{align}
&E\left[\int_0^T\left|\frac{f_s\circ\tau_{\epsilon h}-f_s}{\epsilon}\right|^rds\right]=E\left[\int_0^T\left|\frac{1}{\epsilon}\int_0^\epsilon <\dm f_s\circ\tau_{u h} ,h>_{H}du\right|^rds\right]\notag\\
&\le \frac{\epsilon^{r-1}}{\epsilon^r}\int_0^T\int_0^\epsilon E\left[\left| <\dm f_s ,h>_{H}\right|^r\circ\tau_{u h}\right] duds=\frac{1}{\epsilon}\int_0^T\int_0^\epsilon E\left[\left|<\dm f_s ,h>_{H}\right|^r\mathcal{E}(u\dot{h})(T)\right] duds\notag\\
&\le \frac{1}{\epsilon}\int_0^T\left(E\left[\left|<\dm f_s ,h>_{H}\right|^q\right]\right)^{r/q}\int_0^\epsilon \left(E\left[\left|\mathcal{E}(u\dot{h})(T)\right|^{q/(q-r)}\right]\right)^{(q-r)/q} duds\notag\\
&\le \underset{u\in (0,1)}{\sup} \left(E\left[\left|\mathcal{E}(u\dot{h})(T)\right|^{q/(q-r)}\right]\right)^{(q-r)/q}\int_0^T\left(E\left[\left|<\dm f_s ,h>_{H}\right|^q\right]\right)^{r/q} ds\label{equi3}
\end{align}
From Jensen  inequality  we get
 \begin{align}
&\int_0^T\left(E\left[\left|<\dm f_s ,h>_{H}\right|^q\right]\right)^{r/q} ds\le \left(\|h\|_{H}^{r}\int_0^TE\left[\|\dm f_s\|_H^q\right]ds\right)^{r/q}<\infty\label{equi4}
\end{align}
Using Proposition \ref{cama} and replacing \eqref{equi4} and \eqref{equi3} in \eqref{equi1} we obtain
\begin{align}
&\underset{\epsilon\in(0,1)}{\sup}E\left[\left\|\frac{f\circ\tau_{\epsilon h}-f}{\epsilon}\right\|_{L^p([0,T])}^{p+\eta}\right]<\infty\label{equi}
\end{align}

The convergence in probability \eqref{pr1} and \eqref{equi} imply \cite[Theorem 6.6.2]{Resnik:1999}
\begin{align*}
\frac{f\circ\tau_{\epsilon h}-f}{\epsilon}~\overset{L^p(\Omega, L^p([0,T])}{\underset{\epsilon\rightarrow 0}{\rightarrow}}F[h]
\end{align*}
This implies
\begin{align*}
\underset{\epsilon\rightarrow 0}{\lim}\int_0^TE\left[\left|\frac{f_s\circ\tau_{\epsilon h}-f_s}{\epsilon}-F_s[h]\right|^p \right]ds=0
\end{align*}
Since convergence in $L^1([0,T])$ implies convergence in probability, there exist an sequence $\{\epsilon_n\}_{n\in \N}$ such that $\epsilon_n\rightarrow 0$ when $n\rightarrow \infty$ and for almost all $s\in[0,T]$
\begin{align*}
\underset{\epsilon_n\rightarrow 0}{\lim}E\left[\left|\frac{f_s\circ\tau_{\epsilon_n h}-f_s}{\epsilon_n}-F_s[h]\right|^p \right]=0
\end{align*}
But from this and \eqref{eqsa} we get for almost all $s\in[0,T]$
\begin{align*}
\frac{f_s\circ\tau_{\epsilon_n h}-f_s}{\epsilon_n}\underset{\epsilon_n\rightarrow 0}{\overset{L^p(\Omega)}{\rightarrow}}F_s[h], ~\frac{f_s\circ\tau_{\epsilon_n h}-f_s}{\epsilon_n}\underset{\epsilon_n\rightarrow 0}{\overset{L^p(\Omega)}{\rightarrow}}<\dm f_s,h>_{H}
\end{align*}
Thus $F_s(\omega)[h]=<\dm f_s(\omega),h(\omega)>_{H}$ for almost all  $s\in[0,T]$ and $\omega\in \Omega$
\end{proof}

\section*{Appendix B: The proof of lemma \ref{lemaA22}}
\label{apendixb}
\begin{proof}
Inequalities \eqref{eql1} imply that $f\in \D ^{1,q}(C([0,T])$. Thus, from Theorem \ref{thg1} we  get that $f$ is Ray Absolutely Continuous and Stochastically Gate{\^a}ux differentiable. Thus $f\in L^q(\Omega,C([0,T])) $ and there exists $F\in L^q(\Omega,L(H,C([0,T]))$ such that
\begin{align*}
\underset{t\in[0,T]}{\sup}\left|\frac{f_t\circ\tau_{\epsilon h}-f_t}{\epsilon}-F_t[h]\right|~\overset{\p}{\underset{\epsilon\rightarrow 0}{\rightarrow}}~0
\end{align*}
This implies that for any fixed $t\in[0,T]$ we have
\begin{align}
\left|\frac{f_t\circ\tau_{\epsilon h}-f_t}{\epsilon}-F_t[h]\right|~\overset{\p}{\underset{\epsilon\rightarrow 0}{\rightarrow}}~0\label{cop1}
\end{align}
But for any fixed $t\in[0,T]$, from inequalities \eqref{eql1} we have $f_t\in \D ^{1,q}$, so from Theorem \ref{thg1} we  get that $f_t$ is Ray Absolutely Continuous and Stochastically Gate{\^a}ux differentiable.  From \cite{MastroliaPossamaiRveillac:2017} we have
\begin{align}
\left|\frac{f_t\circ\tau_{\epsilon h}-f_t}{\epsilon}-<\dm f_t, g>_H\right|~\overset{\p}{\underset{\epsilon\rightarrow 0}{\rightarrow}}~0\label{cop2}
\end{align}
Proceeding as at the end of the proof of Lemma \ref{lemaA21} and taking a sequence $\{\epsilon_n\}_{n\in \N}$ such that $\epsilon_n\rightarrow 0$ when $n\rightarrow \infty$ from the two convergence in probability \eqref{cop1} and \eqref{cop2} we get that $\frac{f_t\circ\tau_{\epsilon_n h}-f_t}{\epsilon_n}$ converges almost surely to both $<\dm f_t, g>_H$ and $F_t[h]$, so for any $t\in[0,T]$, $F_t[h]=<\dm f_t, g>_H$, $\p$ a.s. and we have \eqref{probi11}.

Proceeding as in the proof of Lemma 4.2 in \cite{MastroliaPossamaiRveillac:2017}, using \eqref{eqsh},  H\"{o}lder  inequality, the Cameron-Martin Formula, and inequalities \eqref{eql1}   we obtain
\begin{align*}
&\underset{\epsilon\in(0,1]}{\sup}E\left[\underset{t\in[0,T]}{\sup}\left|\frac{f_t\circ\tau_{\epsilon h}-f_t}{\epsilon}\right|^2\right]=\underset{\epsilon\in(0,1]}{\sup}E\left[\underset{t\in[0,T]}{\sup}\left|\frac{1}{\epsilon}\int_0^\epsilon <\dm f_t\circ\tau_{uh},h>_Hdu\right|^2\right]\\
&\le\underset{\epsilon\in(0,1]}{\sup}\frac{1}{\epsilon}\int_0^\epsilon E\left[\underset{t\in[0,T]}{\sup} \left|<\dm f_t,h>_H\right|^2\circ\tau_{uh}\right]du=\underset{\epsilon\in(0,1]}{\sup}\frac{1}{\epsilon}\int_0^\epsilon E\left[\underset{t\in[0,T]}{\sup} \left|<\dm f_t,h>_H\right|^2\mathcal{E}(u\dot{h})(T)\right] du\notag\\
&\le\underset{u\in(0,1)}{\sup}\left(E\left[\left(\mathcal{E}(u\dot{h})(T)\right)^{q/(q-2)}\right]\right)^{(q-2)/q}  \left(E\left[\left(\underset{t\in[0,T]}{\sup} \left|<\dm f_t,h>_H\right|^2\right)^{q/2}\right]\right)^{2/q} \notag\\
&=\|h\|_H^2\underset{u\in(0,1)}{\sup}\left(E\left[\left(\mathcal{E}(u\dot{h})(T)\right)^{q/(q-2)}\right]\right)^{(q-2)/q}  \left(E\left[\left(\underset{t\in[0,T]}{\sup} \int_0^T|\dm_s f_t|^2ds\right)^{q/2}\right]\right)^{2/q}<\infty 
\end{align*}
\end{proof}
\section*{Appendix C: The proof of lemma \ref{lemh1}}
\label{aplemma3}
\begin{proof}
The proof is similar  with the proof of Theorem 3.21 in \cite{ImkellerReisSalkeld:2019}. Equation \eqref{eq4bis} is the same with equation  (3.6) in  \cite{ImkellerReisSalkeld:2019} with
\begin{align*}
&A(s,\om):=\dm_s\alpha (r,\om)\mathds{1}_{\{r\ge s\}}+ \left(\nabla_x\sigma(X(s)(\om))Y(s)(\om)+V(s,\om)\right)\mathds{1}_{\{r\le s\}}\text{ instead of } \sigma(X(s)(\om)) \text{ in (3.6) in  \cite{ImkellerReisSalkeld:2019}}\\
&F(s,u,\omega):=\dm _s U(u,\omega)+\dm_s\nabla_x b(X(u))Y(u) \text{ instead of } U(s,u,\omega)\text{ in (3.6) in  \cite{ImkellerReisSalkeld:2019}}\\
&G(s,u,\omega):=\dm _s V(u,\omega)+\dm _s\nabla_x\sigma(X(u))Y(u)\text{ instead of } V(s,u,\omega)\text{ in (3.6) in  \cite{ImkellerReisSalkeld:2019}}.
\end{align*}
Notice that  $F:[0,T]^2\times \Omega\rightarrow \R^{d\times m}$, $G:[0,T]^2\times \Omega\rightarrow \R^{(d\times m)\times m}$, and $F(s,u,\omega)=0$ and $G(s,u,\omega)=0$ for $s>u$.

Applying twice  Minkovski  inequality, H\"{o}lder inequality, assumption {\bf{C}}, \eqref{alp2}, \eqref{V1}, and \eqref{mom3y} we get $E\left[\|A(\cdot)\|_2^p\right]\le \infty$.

We have assumed that for any $t\in [r,T]$, $U(t)\in \D^{1,\infty}(\R^d)$ and using  \eqref{U1}, \eqref{U2} we obtain
\begin{align}
&E\left[\left(\int_r^T\left(\int_0^T|\dm_ s U(u)|^2ds\right)^{1/2}du\right)^p\right]\le C E\left[\underset{s\in [r,T]}{\sup}\left(\int_0^T|\dm_t U(s)|^2dt\right)^{p/2} \right]<\infty\label{con4U}
\end{align}

From assumption {\bf{C}}, \eqref{mom3}, \eqref{malmom1},\eqref{us2}, and \eqref{mom3y} we get 
\begin{align}
&E\left[\left(\int_r^T\left(\int_0^T|\dm_ s \nabla_x b(X(u))Y(u)|^2ds\right)^{1/2}du\right)^p\right]<\infty\label{con4B}
\end{align}
Using twice  Minkovski  inequality, \eqref{con4U}, and \eqref{con4B} we obtain $E\left[\left(\int_r^T\|F(\cdot,u)\|_2du\right)^p\right]<\infty$.

Similarly, from \eqref{V1}, \eqref{V2} we get
\begin{align}
&E\left[\left(\int_r^T\int_0^T|\dm_ s V^j(u,\omega)|^2ds du\right)^{p/2}\right]&\le C E\left[\underset{u\in[r,T]}{\sup}\left(\int_0^T|\dm_ s V(u,\omega)|^2ds\right)^{p/2} \right]<\infty.\label{con4V}
\end{align}
Using twice  Minkovski  inequality, \eqref{con4V}, \eqref{con4s} we obtain $E\left[\left(\int_r^T\|G(\cdot,u)\|_2^2du\right)^{p/2}\right]<\infty$.

From Theorem 2.5 in \cite{ImkellerReisSalkeld:2019} we obtain
\begin{align}
&E\left[\left(\underset{t\in [r\vee s,T]}{\sup}\int_0^T|M_s(t)|^2ds\right)^{p/2}\right]\le E\left[\left(\underset{t\in [r\vee s,T]}{\sup}\|M_\cdot(t)\|_2\right)^{p}\right]\notag\\
&\le E\left[\|A(\cdot)\|_2^p\right]+E\left[\left(\int_r^T\|F(\cdot,u)\|_2du\right)^p\right]+E\left[\left(\int_r^T\|G(\cdot,u)\|_2^2du\right)^{p/2}\right]<\infty\notag
\end{align}
\end{proof}
\section*{Appendix D: The proof of lemma \ref{lemah2}}
\label{aplemma4}
\begin{proof}
Notice that from assumption {\bf{C}} and \eqref{V1} we get $V$, $\Sigma^i\in H_m^2$, so using \eqref{helpto} we have $\mathbb{P}$ a.s.
\begin{align}
& \frac{1}{\ep}\left(\int_r^TV^i(u,\omega+\ep h)dW^i(u)(\om+\ep h)-\int_r^TV^i(u,\omega)dW^i(u)(\om)\right)\notag\\
&=\int_r^TV^i(u,\omega+\ep h) \dot{h}^i(u)du+\int_r^T\frac{V^i(u,\omega+\ep h)-V^i(u,\omega)}{\ep}dW^i(u)(\om),\label{eqvop}\\
&\frac{1}{\ep}\biggl(\int_r^T\Sigma^i(u,\om +\ep h)Y(u)(\om +\ep h)dW^i(u)(\om +\ep h)-\int_r^T\Sigma^i(u,\om)Y(u)(\om)dW^i(u)(\om)\biggl)\notag\\
&=\int_r^T\Sigma^i(u,\om +\ep h)Y(u)(\om +\ep h)\dot{h}^i(u)du\notag\\
&+ \int_r^T\frac{\Sigma^i(u,\om +\ep h)Y(u)(\om +\ep h)-\Sigma^i(u,\om)Y(u)(\om)}{\ep}dW^i(u)(\om).\label{eqsop}
\end{align}
From  \eqref{helpto}, \eqref{eqvop}, and \eqref{eqsop} we get
\begin{align*}
&P_{\ep}(t)(\om)=\frac{\alpha (r,\om+\ep h)-\alpha (r,\om)}{\ep}+\int_r^t \frac{U(u,\omega+\ep h)- U(u,\omega)}{\ep}du\\
&+\sum_{i=1}^m\int_r^tV^i(u,\omega+\ep h) \dot{h}^i(u)du+\sum_{i=1}^m\int_r^t\Sigma^i(u,\om +\ep h)Y(u)(\om +\ep h)\dot{h}^i(u)du\\
&+\int_r^t\frac{\left(B(u,\om +\ep h))-B(u,\om))\right)Y(u)(\om )}{\ep}du+\int_r^t\frac{B(u,\om +\ep h))\left(Y(u)(\om+\ep h )-Y(u)(\om)\right)}{\ep}du\\
&+\sum_{i=1}^m\int_r^t\frac{V^i(u,\omega+\ep h)-V^i(u,\omega)}{\ep}dW^i(u)(\om)+ \sum_{i=1}^m\int_r^t\frac{\left(\Sigma^i(u,\om +\ep h)-\Sigma^i(u,\om)\right)Y(u)(\om)}{\ep}dW^i(u)(\om)\\
&+\sum_{i=1}^m\int_r^t\frac{\Sigma^i(u,\om +\ep h)\left(Y(u)(\om+\ep h)-Y(u)(\om)\right)}{\ep}dW^i(u)(\om)
\end{align*} 
Using It{\^o}'s formula for $f(x)=x^2$ we have
\begin{align}
&|P_{\ep}(t)(\om)|^2=\left|\frac{\alpha (r,\om+\ep h)-\alpha (r,\om)}{\ep}\right|^2\label{e1}\\
&+2\int_r^t\left\langle P_{\ep}(u)(\om),\frac{U(u,\omega+\ep h)- U(u,\omega)}{\ep}\right\rangle du\label{e2}\\
&+2\sum_{i=1}^m\int_0 ^t\left\langle P_{\ep}(u)(\om),V^i(u,\omega+\ep h) \dot{h}^i(u)\right\rangle du\label{e3}\\
&+2 \sum_{i=1}^m\int_r^t\left\langle P_{\ep}(u)(\om),\Sigma^i(u,\om +\ep h)Y(u)(\om +\ep h)\dot{h}^i(u)\right\rangle du\label{e4}\\
&+2\int_r^t\left\langle P_{\ep}(u)(\om),\frac{\left(B(u,\om +\ep h))-B(u,\om))\right)Y(u)(\om )}{\ep}\right\rangle du\label{e5}\\
&+2\int_r^t\left\langle P_{\ep}(u)(\om),\frac{B(u,\om +\ep h))\left(Y(u)(\om+\ep h )-Y(u)(\om)\right)}{\ep}\right\rangle du\label{e6}\\
&+2\sum_{i=1}^m\int_r^t\left\langle P_{\ep}(u)(\om),\frac{V^i(u,\omega+\ep h)-V^i(u,\omega)}{\ep}\right\rangle dW^i(u)(\om)\label{e7}\\
&+2\sum_{i=1}^m\int_r^t\left\langle P_{\ep}(u)(\om),\frac{\left(\Sigma^i(u,\om +\ep h)-\Sigma^i(u,\om)\right)Y(u)(\om)}{\ep}\right\rangle dW^i(u)(\om)\label{e8}\\
&+2\sum_{i=1}^m\int_r^t\left\langle P_{\ep}(u)(\om),\frac{\Sigma^i(u,\om +\ep h)\left(Y(u)(\om+\ep h)-Y(u)(\om)\right)}{\ep}\right\rangle dW^i(u)(\om)\label{e9}\\
&+\sum_{i=1}^m\int_r^t\biggl|\frac{V^i(u,\omega+\ep h)-V^i(u,\omega)}{\ep}+\frac{\left(\Sigma^i(u,\om +\ep h)-\Sigma^i(u,\om)\right)Y(u)(\om)}{\ep}\notag\\
&+ \frac{\Sigma^i(u,\om +\ep h)\left(Y(u)(\om+\ep h)-Y(u)(\om)\right)}{\ep}\biggl| ^2du\label{e10}
\end{align}
We take a supremum over $t$ then expectations. Let $n$ be a positive integer that we will choose later. We find the following upper bounds.

For \eqref{e1} we have from Cauchy-Schwartz inequality, \eqref{con1a}, and \eqref{alp2} with $p=2$ 
\begin{align*}
&E\left[\left|\frac{\alpha (r,\om+\ep h)-\alpha (r,\om)}{\ep}\right|^2\right]\le 2E\left[\left|\frac{\alpha (r,\om+\ep h)-\alpha (r,\om)}{\ep}-\dm \alpha (r) [h]\right|^2\right]+2E\left[\left|\dm \alpha (r) [h]\right|^2\right]\\
&\le o(1)+2\|\dot{h}\|_2^2E\left[\underset{r \in [0,T]}{\sup}\int_0^T|\dm_s \alpha (r)|^2ds\right]\text{ as }\ep\rightarrow 0.
\end{align*}
For \eqref{e2} we have from Cauchy-Schwartz, Jensen, and Young  inequalities, and from \eqref{con3U}
\begin{align*}
&E\left[\underset{t \in [r,T]}{\sup}\left|2\int_r^t\left\langle P_{\ep}(u)(\om),\frac{U(u,\omega+\ep h)- U(u,\omega)}{\ep}\right\rangle du\right|\right]\\
&\le\frac{1}{n}E\left[\|P_{\ep}\|_\infty^2\right]+nT\underset{\epsilon\in(0,1]}{\sup}E\left[\underset{u \in [r,T]}{\sup}\left|\frac{U(u,\omega+\ep h)- U(u,\omega)}{\ep}\right|^2\right]
\end{align*}
To find upper bounds for \eqref{e5} we first notice that  from \eqref{con6B}, using also Cauchy-Schwartz and Young  inequalities, we get 
\begin{align*}
&E\left[\underset{t \in [r,T]}{\sup}\left|2\int_r^t\left\langle P_{\ep}(u)(\om),\frac{\left(B(u,\om +\ep h))-B(u,\om))\right)Y(u)(\om )}{\ep}\right\rangle du\right|\right]\\
&=o(1)+E\left[\frac{\|P_{\ep}\|_\infty^2}{n}\right]+n2d^2T\|\dot{h}\|_2^2\sum_{i,j=1}^d \left(E\left[\underset{u\in[r,T]}{\sup}|Y(u)(\om )|^4du\right]\right)^{1/2}\\
&\left(E\left[\int_r^T\left(\int_0^T|\dm_s B^{i,j}(u)(\om )|^2ds\right)^2du \right]\right)^{1/2}
\end{align*}
For \eqref{e6} we have from \eqref{c4}
\begin{align*}
&2E\left[\underset{t \in [r,T]}{\sup}\left|\int_r^t\left\langle P_{\ep}(u)(\om),\frac{B(u,\om +\ep h))\left(Y(u)(\om+\ep h )-Y(u)(\om)\right)}{\ep}\right\rangle du\right|\right]\\
&= 2E\left[\underset{t \in [r,T]}{\sup}\int_r^t\left\langle P_{\ep}(u)(\om),B(u,\om +\ep h))P_{\ep}(u)\right\rangle du\right]\le2L\int_r^TE\left[ \|P_{\ep}(s)(\om)\|^2_{\infty,u}\right]du
\end{align*}
For \eqref{e7} we have from \eqref{con3V0}, and Burkholder-Davis-Gundy, Cauchy-Schwartz, and Young inequalities  
\begin{align*}
&2E\left[\underset{t \in [r,T]}{\sup}\left|\int_r^t\left\langle P_{\ep}(u)(\om),\frac{V^i(u,\omega+\ep h)-V^i(u,\omega)}{\ep}\right\rangle dW^i(u)(\om)\right|\right]\\
&\le 2CE\left[\left(\int_r^T\left|\left\langle P_{\ep}(u)(\om),\frac{V^i(u,\omega+\ep h)-V^i(u,\omega)}{\ep}\right\rangle\right|^2du\right)^{1/2}\right]\\
&\le\frac{1}{n}E\left[\| P_{\ep}\|_{\infty}^2\right]+C^2nT\underset{\epsilon\in(0,1]}{\sup}E\left[\underset{u \in [r,T]}{\sup}\left|\frac{V^i(u,\omega+\ep h)-V^i(u,\omega)}{\ep}\right|^2\right]
\end{align*}
Similarly with the work done for \eqref{e7}, for \eqref{e8} and \eqref{e9} we have from \eqref{con6s0}, assumption {\bf{C}}, Burkholder-Davis-Gundy, Cauchy-Schwartz, and Young  inequalities and 
\begin{align*}
&2E\left[\underset{t \in [r,T]}{\sup}\left|\int_r^t\left\langle P_{\ep}(u)(\om),\frac{\left(\Sigma^i(u,\om +\ep h)-\Sigma^i(u,\om)\right)Y(u)(\om)}{\ep}\right\rangle dW^i(u)(\om)\right|\right]\\
&\le o(1)+\frac{1}{n}E\left[\| P_{\ep}\|_{\infty}^2\right]+2C^2n\|\dot{h}\|_2^2\sum_{j,k=1}^d\left(E\left[\int_r^T\left(\int_r^T\left|\dm_s \Sigma^{i,j,k}(u)\right|^2ds\right)^2du\right]\right)^{1/2}\\
&\left(E\left[\underset{u\in[0,T]}{\sup}|Y(u)|^4du\right]\right)^{1/2},
\end{align*}
\begin{align*}
&2E\left[\underset{t \in [r,T]}{\sup}\left|\int_r^t\left\langle P_{\ep}(u)(\om),\frac{\Sigma^i(u,\om +\ep h)\left(Y(u)(\om+\ep h)-Y(u)(\om)\right)}{\ep}\right\rangle dW^i(u)(\om)\right|\right]\\
&=2E\left[\underset{t \in [r,T]}{\sup}\left|\int_r^t\left\langle P_{\ep}(u)(\om),\Sigma^i(u,\om +\ep h)P_{\ep}(u)(\om)\right\rangle dW^i(u)(\om)\right|\right]\\
&\le C_1\int_r^TE\left[ \underset{s \in [r,u]}{\sup}|P_{\ep}(s)(\om)|^2\right]du+C_1=C_1\int_r^TE\left[ \|P_{\ep}|_{\infty,s}^2\right]du+C_1
\end{align*}
For \eqref{e10} from assumption {\bf{C}}, \eqref{con3V0}, and \eqref{con6s0}
\begin{align*}
&E\biggl[\underset{t \in [r,T]}{\sup}\int_r^t\biggl|\frac{V^i(u,\omega+\ep h)-V^i(u,\omega)}{\ep}+\frac{\left(\Sigma^i(u,\om +\ep h)-\Sigma^i(u,\om)\right)Y(u)(\om)}{\ep}\\
&+ \frac{\Sigma^i(u,\om +\ep h)\left(Y(u)(\om+\ep h)-Y(u)(\om)\right)}{\ep}\biggl| ^2du\biggl]\\
&\le o(1)+ 3T\underset{\epsilon\in(0,1]}{\sup}E\left[\underset{u \in [r,T]}{\sup}\left|\frac{V^i(u,\omega+\ep h)-V^i(u,\omega)}{\ep}\right|^2\right]\\
&+6\|\dot{h}\|_2^2\sum_{j,k=1}^d\left(E\left[\int_r^T\left(\int_0^T\left|\dm_s \Sigma^{i,j,k}(u)\right|^2ds\right)^2du\right]\right)^{1/2}\left(E\left[\underset{u\in[r,T]}{\sup}|Y(u)|^4du\right]\right)^{1/2}\\
&+3C\int_r^TE\left[ \|P_{\ep}\|_{\infty,s}^2\right]du
\end{align*}
For \eqref{e3}, for any $i=1,\ldots, m$ we have from \eqref{con3V0}, Cauchy-Schwartz, and Young  inequalities 
\begin{align*}
&2E\left[\underset{t \in [r,T]}{\sup}\left|\int_r ^t\left\langle P_{\ep}(u)(\om),V^i(u,\omega+\ep h) \dot{h}^i(u)\right\rangle du\right|\right]\\
&\le 2\ep E\left[\int_r^T| P_{\ep}(u)(\om)|\left|\frac{V^i(u,\omega+\ep h)-V^i(u,\omega)}{\ep}\right| |\dot{h}^i(u)| du\right]\\
&\le \frac{2}{n}E\left[\| P_{\ep}\|_\infty^2\right]+n\ep ^2T\|\dot{h}^i\|_2^2\underset{\ep \in (0,1]}{\sup}E\left[\underset{u \in [r,T]}{\sup} \left|\frac{V^i(u,\omega+\ep h)-V^i(u,\omega)}{\ep}\right|^2 \right]\\
&+n\| \dot{h}^i\|^2_2E\left[\underset{u \in [r,T]}{\sup} |V^i(u,\omega)|^2 du\right]
\end{align*}
Proceeding as for \eqref{e3} and using assumption {\bf{C}}, \eqref{con6s0}, \eqref{mom3y} with $q=2$ and 4, for \eqref{e4} we have
\begin{align*}
&2E\left[\underset{t \in [r,T]}{\sup}\left|\int_r^t\left\langle P_{\ep}(u)(\om),\Sigma^i(u,\om +\ep h)Y(u)(\om +\ep h)\dot{h}^i(u)\right\rangle du\right|\right]\\
&\le \ep E\left[\int_r^T\left|\left\langle P_{\ep}(u)(\om),\frac{\Sigma^i(u,\om +\ep h)-\Sigma^i(u,\om)}{\ep}Y(u)(\om)\dot{h}^i(u)\right\rangle\right| du\right]\\
&+\ep E\left[\int_r^T\left|\left\langle P_{\ep}(u)(\om),\frac{\Sigma^i(u,\om +\ep h)\left(Y(u)(\om +\ep h)-Y(u)(\om)\right)}{\ep}\dot{h}^i(u)\right\rangle\right| du\right]\\
&+E\left[\int_r^T\left|\left\langle P_{\ep}(u)(\om),\Sigma^i(u,\om )Y(u)(\om )\dot{h}^i(u)\right\rangle\right| du\right]\\
&\le\frac{2}{n}E\left[\| P_{\ep}\|_\infty^2\right]+\ep CE\left[\int_r^T\left\|P_{\ep}\right\|_{\infty,u}^2|\dot{h}^i(u)| du\right]+o(1)+nC\| \dot{h}^i\|^2_2E\left[\underset{u \in [r,T]}{\sup}|Y(u)(\om )|^2du\right]\\
&+2n\ep ^2\|\dot{h}^i\|_2^4
\sum_{j,k=1}^d\left(E\left[\int_r^T\left(\int_r^T\left|\dm_s \Sigma^{i,j,k}(u)\right|^2ds\right)^2du\right]\right)^{1/2}\left(E\left[\underset{u\in[r,T]}{\sup}|Y(u)|^4du\right]\right)^{1/2}
\end{align*}

Combining all these inequalities and choosing $n=9$ we have
\begin{equation*}
\frac{1}{9}E\left[\| P_{\ep}\|_\infty^2\right]\le E\left[\|A_\ep\|^2_\infty\right]+C_1\int_r^T(1+\ep|h^i(u)|)E\left[\left\|P_{\ep}\right\|_{\infty,u}^2\right]du,
\end{equation*}
where $E\left[\|A_\ep\|^2_\infty\right]=O(1)$ as $\ep \rightarrow 0$. Gr{\"o}nwall's inequality yields that $E\left[\| P_{\ep}\|_\infty^2\right]=O(1)$ as $\ep \rightarrow 0$.
\end{proof}
\section*{Appendix E: The proof of lemma \ref{lemh3}}
\label{aplemma5}
\begin{proof}
Let $t\in [0,T]$. We introduce the notation:
\begin{align*}
&Z_\ep(t)(\om):=P_{\ep}(t)(\om)-M^h(t)(\om)
\end{align*}
We write the SDE for $Z_\ep(t)(\om)$:
\begin{align}
&Z_\ep(t)(\om)=\frac{\alpha (r,\om+\ep h)-\alpha (r,\om)}{\ep}-\int_0^r\dm_s\alpha (r,\om)\dot{h}(s)ds\label{pe1}\\
&+\int_r^t \left(\frac{U(u,\omega+\ep h)- U(u,\omega)}{\ep}-\int_0^u \dm _s U(u,\omega)\dot{h}(s)ds\right)du\label{pe2}\\
&+\ep\sum_{i=1}^m\int_r^t\frac{V^i(u,\omega+\ep h) -V^i(u,\om)}{\ep}\dot{h}^i(u)du\label{pe3}\\
&+\ep\sum_{i=1}^m\int_r^t\frac{\Sigma^i(u,\om +\ep h)-\Sigma^i(u,\om)}{\ep} Y(u)(\om)\dot{h}^i(u)du\label{pe4}\\
&+\ep\sum_{i=1}^m\int_r^t\Sigma^i(u,\om +\ep h)\frac{Y(u)(\om +\ep h)-Y(u)(\om)}{\ep} \dot{h}^i(u)du\label{pe5}\\
&+\int_r^t\left(\frac{\left(B(u,\om +\ep h))-B(u,\om))\right)Y(u)(\om )}{\ep}-\int_0^u\dm_s \nabla_x b(X(u))Y(u)\dot{h}(s)ds\right)du\label{pe6}\\
&+\sum_{i=1}^m\int_r^t\left(\frac{V^i(u,\omega+\ep h)-V^i(u,\omega)}{\ep}-\int_0^u\dm _s V^i(u,\omega)\dot{h}(s)ds\right)dW^i(u)(\om)\label{pe7}\\
&+\sum_{i=1}^m\int_r^t\left(\frac{\left(\Sigma^i(u,\om +\ep h)-\Sigma^i(u,\om)\right)Y(u)(\om)}{\ep}-\int_0^u\dm _s \nabla_x\sigma^i(X(u))Y(u)\dot{h}(s)ds\right)dW^i(u)(\om)\label{pe8}\\
&+\int_r^t B(u,\om))\left(\frac{\left(Y(u)(\om+\ep h )-Y(u)(\om)\right)}{\ep}-\int_0^u M_s(u)\dot{h}(s)ds\right)du\label{pe9}\\
&+\ep\int_r^t\frac{\left(B(u,\om +\ep h))-B(u,\om))\right)P_\ep(u)(\om )}{\ep}du\label{pe10}\\
&+\sum_{i=1}^m\int_r^t\Sigma^i(u,\om)\left(\frac{Y(u)(\om+\ep h)-Y(u)(\om)}{\ep}-\int_0^u M_s(u)\dot{h}(s)ds\right)dW^i(u)(\om)\label{pe11}\\
&+\ep\sum_{i=1}^m\int_r^t\frac{\left(\Sigma^i(u,\om +\ep h)-\Sigma^i(u,\om)\right)P_\ep(u)(\om)}{\ep}dW^i(u)(\om)\label{pe12}
\end{align}
We take supremum over $t\in [r,T]$, and next we prove mean  convergence for each  term.

For \eqref{pe10}  from \eqref{con3B} and Lemma \ref{lemah2}, using Cauchy-Schwartz inequality we have
\begin{align*}
&E\left[\underset{t \in [r,T]}{\sup}\left|\ep\int_r^t\frac{\left(B(u,\om +\ep h))-B(u,\om))\right)P_\ep(u)(\om )}{\ep}du\right|\right]\\
&\le \ep C\left(E\left[\underset{u \in [r,T]}{\sup}\left|\frac{B(u,\om +\ep h))-B(u,\om))}{\ep}\right|^2du\right]\right)^{1/2}\left(E\left[\underset{u \in [r,T]}{\sup}| P_\ep(u)(\om )|^2du\right]\right)^{1/2}\le O(\ep)
\end{align*}
For \eqref{pe12} for any $i=1,\ldots,m$ we have from Burkholder-Davis-Gundy inequality, \eqref{con3s0}, and Lemma \ref{lemah2}, using Cauchy-Schwartz inequality we get
\begin{align*}
&\ep E\left[\underset{t \in [r,T]}{\sup}\left|\int_r^t\frac{\left(\Sigma^i(u,\om +\ep h)-\Sigma^i(u,\om)\right)P_\ep(u)(\om)}{\ep}dW^i(u)(\om)\right|\right]\\
&\le C\ep E\left[\left(\int_r^T\left|\frac{\Sigma^i(u,\om +\ep h)-\Sigma^i(u,\om)}{\ep}\right|^2|P_\ep(u)(\om)|^2du\right)^{1/2}\right]\\
&\le C\ep \left(E\left[\underset{u \in [r,T]}{\sup}\left|\frac{\Sigma^i(u,\om +\ep h)-\Sigma^i(u,\om)}{\ep}\right|^2\right]\right)^{1/2}\left(E\left[\underset{u \in [r,T]}{\sup}|P_\ep(u)(\om)|^2du\right]\right)^{1/2}\le O(\ep)
\end{align*}
For \eqref{pe1} we have from \eqref{con1a}
\begin{align*}
\underset{\ep\rightarrow\infty}{\lim}E\left[\left|\frac{\alpha (r,\om+\ep h)-\alpha (r,\om)}{\ep}-\int_0^r\dm_s\alpha (r,\om)\dot{h}(s)ds\right|\right]=0
\end{align*}
For \eqref{pe2} we have  from  \eqref{con1U} with $p=2$
\begin{align*}
&0\le E\left[\left(\underset{t \in [r,T]}{\sup}\left|\int_r^t \left(\frac{U(u,\omega+\ep h)- U(u,\omega)}{\ep}-\int_0^u \dm _s U(u,\omega)\dot{h}(s)ds\right)du\right|\right)^2\right]\\
&\le E\left[\left(\int_r^T \left|\frac{U(u,\omega+\ep h)- U(u,\omega)}{\ep}-\int_0^u \dm _s U(u,\omega)\dot{h}(s)ds\right|du\right)^2\right]\rightarrow 0 \text{ when } \ep\rightarrow 0.
\end{align*}
Similarly from \eqref{con5B} we get for \eqref{pe6}
\begin{align*}
&\underset{\ep \rightarrow 0}{\lim}E\biggl[\biggl(\underset{t \in [r,T]}{\sup}\biggl|\int_r^t\biggl(\frac{\left(B(u,\om +\ep h))-B(u,\om))\right)Y(u)(\om )}{\ep}-\int_0^u\dm_s \nabla_x b(X(u))Y(u)\dot{h}(s)ds\biggl)du\biggl|\biggl)^2\biggl]=0
\end{align*}
For \eqref{pe7} for any $i=1,\ldots,m$, from \eqref{con1V0} with $p=2$, Burkholder-Davis-Gundy inequality,  and  using Cauchy-Schwartz inequality we have
\begin{align*}
 &0\le E\left[\underset{t \in [r,T]}{\sup}\left|\int_r^t\left(\frac{V^i(u,\omega+\ep h)-V^i(u,\omega)}{\ep}-\int_0^u\dm _s V^i(u,\omega)\dot{h}(s)ds\right)dW^i(u)(\om)\right|\right]\\
&\le C\left( E\left[\int_r^T\left|\frac{V^i(u,\omega+\ep h)-V^i(u,\omega)}{\ep}-\int_0^u\dm _s V^i(u,\omega)\dot{h}(s)ds\right|^2d(u)\right]\right)^{1/2}\underset{\ep \rightarrow 0}{\rightarrow}0
\end{align*}
Thus
\begin{align*}
\underset{\ep \rightarrow 0}{\lim}  E\left[\left\|\int_r^t\left(\frac{V^i(u,\omega+\ep h)-V^i(u,\omega)}{\ep}-\int_0^u\dm _s V^i(u,\omega)\dot{h}(s)ds\right)dW^i(u)(\om)\right\|_\infty\right]=0
\end{align*}
Similarly from \eqref{con5s0} we get for \eqref{pe8} for any $i=1,\ldots,m$
\begin{align*}
&\underset{\ep \rightarrow 0}{\lim}  E\biggl[\biggl\|\int_r^t\biggl(\frac{\left(\Sigma^i(u,\om +\ep h)-\Sigma^i(u,\om)\right)Y(u)(\om)}{\ep}-\int_0^u\dm _s \nabla_x\sigma^i(X(u))Y(u)\dot{h}(s)ds\biggl)dW^i(u)(\om)\biggl\|_\infty\biggl]=0
\end{align*}
For \eqref{pe3} for any $i=1,\ldots,m$ we get from Cauchy-Schwartz inequality and \eqref{con3V0}
\begin{align*}
&\ep^2E\left[\left(\underset{t \in [r,T]}{\sup}\left|\int_r^t\frac{V^i(u,\omega+\ep h) -V^i(u,\om)}{\ep}\dot{h}^i(u)du\right|\right)^2\right]\\
&\le T\ep^2\|\dot{h}^i\|_2^2E\left[\underset{u \in [r,T]}{\sup}\left|\frac{V^i(u,\omega+\ep h) -V^i(u,\om)}{\ep}\right|^2\right]=O(\ep^2)
\end{align*}
For \eqref{pe4} for any $i=1,\ldots,m$ we get from Cauchy-Schwartz inequality, \eqref{con3s0}, and \eqref{mom3y} with $p=2$
\begin{align*}
&\ep E\left[\underset{t \in [r,T]}{\sup}\left|\int_r^t\frac{\Sigma^i(u,\om +\ep h)-\Sigma^i(u,\om)}{\ep} Y(u)(\om)\dot{h}^i(u)du\right|\right]\\
&\le \ep\|\dot{h}^i\|_2 \left(E\left[\underset{u \in [r,T]}{\sup}\left|\frac{\Sigma^i(u,\om +\ep h)-\Sigma^i(u,\om)}{\ep}\right|^2\right]\right)^{1/2}\left(E\left[\underset{u \in [r,T]}{\sup}| Y(u)(\om)|^2\right]\right)^{1/2}=O(\ep)
\end{align*}
For \eqref{pe5} we have for any $i=1,\ldots,m$ from assumption {\bf{C}}, Cauchy-Schwartz inequality, and Lemma \ref{lemah2}
\begin{align*}
&\ep^2 E\left[\left(\underset{t \in [r,T]}{\sup}\left|\int_r^t\Sigma^i(u,\om +\ep h)\frac{Y(u)(\om +\ep h)-Y(u)(\om)}{\ep} \dot{h}^i(u)du\right|\right)^2\right]\\
&\le \ep^2 C^2 E\left[\int_r^T|P_\ep(\om)|^2du\int_r^T| \dot{h}^i(u)|^2du\right]\le T\ep^2 C^2 \| \dot{h}^i\|_2^2E\left[\underset{t \in [r,T]}{\sup}|P_\ep(\om)|^2\right]=O(\ep^2)
\end{align*}
Rewriting \eqref{pe9} and \eqref{pe11} in terms of $Z_\ep$ and using all the previous inequalities we have
\begin{align*}
Z_\ep(t)(\om)=A_\ep(t)(\om)+\int_r^t B(u,\om)Z_\ep(u)(\om)du+\sum_{i=1}^m\int_r^t\Sigma^i(u,\om)Z_\ep(u)(\om)dW^i(u),
\end{align*}
where the sequence $A_\ep$ is such that $\|A_\ep\|_\infty\overset{\p}{\rightarrow}0$ as $\ep\rightarrow 0$. From Proposition 2.6 in \cite{ImkellerReisSalkeld:2019} we get $\|Z_\ep\|_\infty\overset{\p}{\rightarrow}0$ as $\ep\rightarrow 0$.
\end{proof}

}

\end{document}